\documentclass[preprint]{imsart}

\usepackage{amsmath,amssymb,amsthm,amsfonts,amstext,amsbsy,amscd}

\usepackage{comment}
\usepackage{graphicx}
\usepackage{mathrsfs}
\usepackage{xcolor}
\usepackage[authoryear]{natbib}

\RequirePackage[colorlinks,citecolor=blue,urlcolor=blue]{hyperref}

\numberwithin{equation}{section}

\swapnumbers

\newtheorem{satz}{Satz}[section]

\newtheorem{theorem}[satz]{Theorem}
\newtheorem{proposition}[satz]{Proposition}
\newtheorem{corollary}[satz]{Corollary}
\newtheorem{lemma}[satz]{Lemma}

\newtheorem{definition}[satz]{Definition}
\theoremstyle{remark}

\newtheorem{remark}[satz]{Remark}

\newtheorem{example}[satz]{Example}

\DeclareMathOperator{\E}{{\mathbb E}}

\DeclareMathOperator{\R}{{\mathbb R}}
\DeclareMathOperator{\C}{{\mathbb C}}
\DeclareMathOperator{\Z}{{\mathbb Z}}
\DeclareMathOperator{\N}{{\mathbb N}}

\DeclareMathOperator{\HH}{{\mathbb H}}
\DeclareMathOperator{\T}{{\mathbb T}}
\DeclareMathOperator{\PP}{{\mathbb P}}

\DeclareMathOperator{\spann}{span}
\DeclareMathOperator{\trace}{trace}
\DeclareMathOperator{\supp}{supp}

  \DeclareMathOperator{\rank}{rank}

\DeclareMathOperator{\diag}{diag} 
 \DeclareMathOperator{\Id}{Id}

\DeclareMathOperator{\Var}{Var} \DeclareMathOperator{\Cov}{Cov}

\DeclareMathOperator{\vek}{vec}
\DeclareMathOperator{\BIAS}{BIAS}

\providecommand{\eps}{\varepsilon}
\renewcommand{\phi}{\varphi}
\renewcommand{\theta}{\vartheta}
\renewcommand{\subset}{\subseteq}

\renewcommand{\cdot}{{\scriptstyle \bullet} }
\providecommand{\abs}[1]{\lvert #1 \rvert}
\providecommand{\norm}[1]{\lVert #1 \rVert}
\providecommand{\bnorm}[1]{{\Bigl\lVert #1 \Bigr\rVert}}
\providecommand{\babs}[1]{{\Bigl\lvert #1 \Bigr\rvert}}
\providecommand{\scapro}[2]{\langle #1,#2 \rangle}
\providecommand{\bscapro}[2]{\Big\langle #1,#2 \Big\rangle}
\providecommand{\floor}[1]{\lfloor #1 \rfloor}

\renewcommand{\le}{\leqslant}
\renewcommand{\ge}{\geqslant}

\providecommand{\mr}{\color{blue}}

\allowdisplaybreaks[4]

\renewcommand{\cal}{\mathscr}     

\begin{document}

\begin{frontmatter}
\title{Rank tests for time-varying covariance matrices observed under noise}
\runtitle{Rank tests for time-varying covariance matrices observed under noise}

\begin{aug}
\author[A]{\fnms{Markus} \snm{Rei\ss}\ead[label=e1]{mreiss@math.hu-berlin.de}}
\and
\author[B]{\fnms{Lars} \snm{Winkelmann}\ead[label=e2]{lars.winkelmann@fu-berlin.de}}

\address[A]{Institut f\"ur Mathematik, Humboldt-Universit\"at zu Berlin,
\printead{e1}}

\address[B]{School of Business and Economics, Freie Universit\"at Berlin,
\printead{e2}}
\end{aug}

\begin{abstract}
We consider a $d$-dimensional  continuous martingale $X(t)$ with quadratic variation matrix $\langle X\rangle_t=\int_0^t \Sigma(s)\,ds$ and develop tests for the rank of its spot covariance matrix $\Sigma(t)$, $t\in[0,1]$.
The process $X$ is observed under observational noise, as is standard for microstructure noise models in high-frequency finance.
We test the null hypothesis ${\mathcal H}_0:\rank(\Sigma(t))\le r$ against local alternatives ${\mathcal H}_{1,n}:\lambda_{r+1}(\Sigma(t))\ge v_n$, where $\lambda_{r+1}$ denotes the $(r+1)$st eigenvalue  and $v_n\downarrow 0$ as the sample size $n\to\infty$. We construct test statistics based on  eigenvalues of carefully calibrated localized spectral covariance matrix estimates. Critical values are provided  non-asymptotically as well as asymptotically via maximal eigenvalues of Gaussian orthogonal ensembles. The power analysis establishes asymptotic consistency for a separation rate $v_n\thicksim  (\underline\lambda_r^{-1/(\beta+1)}n^{-\beta/(\beta+1)})\wedge n^{-\beta/(\beta+2)}$, depending on the H\"older-regularity $\beta$ of $\Sigma$ and a possible spectral gap $\underline\lambda_r\ge 0$ under ${\mathcal H}_0$. A lower bound shows the optimality of this rate. We discuss why the rate is much faster than conventional estimation rates. The theory is illustrated by simulations and a real data example with German government bonds of varying maturity.
\end{abstract}



\end{frontmatter}

\section{Introduction}

Estimating and testing the rank of a covariance matrix is of key interest in statistics. For i.i.d. observations the spiked covariance model \citep{johnstone2001} gives a prototypical setting where the number of spikes or factors, when interpreted as a factor model, is the rank of a dominating signal matrix corrupted by background or idiosyncratic noise. In a high-dimensional framework \cite{onatski2014} develop tests on the presence of spikes  and \citet{cai2015} establish optimal rates for rank detection. Time-varying covariances appear naturally in time series analysis, but inferring their rank often demands strong modeling assumptions. In a continuous-time setting, however, such assumptions can be relaxed when the covariance evolves continuously over time, offering a more flexible framework for statistical inference.

Our focus lies on continuous-time processes observed at high frequency. The foundational model we consider is a continuous, $d$-dimensional martingale
\begin{equation}\label{EqX} X(t)=X(0)+\int_0^t\sigma(s)\,dB(s), \quad t\in[0,1],
\end{equation}
with $\sigma\in L^2([0,1];\R^{d\times p})$, a $p$-dimensional Brownian motion $B$ and some initial value $X(0)$. The instantaneous or spot covariance matrix of $X$ at time $t$ is then given by
\[ \Sigma(t):=\sigma(t)\sigma(t)^\top\in\R^{d\times d}.\]
If $\rank(\Sigma(t))\le r$ holds for all $t\in[0,1]$, then the process $X$ can be represented
by an $r$-dimensional Brownian motion $B$ and suitable $\sigma(t)\in\R^{d\times r}$. See \citet{jacod2008} for a general semi-martingale framework. In financial econometrics, the maximal rank $r$ corresponds to the number of latent risk factors driving the dynamics of asset prices $X$.
\citet{jacod2013}, \citet{ait2017}, along with the references therein, provide further insight into the economic implications.

When high-frequency observations $(X(i/n))_{i=0,\ldots,n}$ are directly available, the maximal rank of $\Sigma(t)$ can be tested using the eigenvalues of local empirical covariance matrices.
This has been accomplished in \citet{RW23} with an optimal rank detection rate of order ${\mathcal O}(n^{-1})$ against local alternatives. This fast rate, especially under a nonparametric model, can be explained by the intrinsic heteroskedasticity of covariance matrix estimators and their eigenvalues, which are tested around zero. The law of the test statistics, however, degenerates asymptotically, being dominated by a bias term.

Models for inference on the volatility often postulate that high-frequency observations are corrupted by independent noise \citep{asjj2014}. A  standard model is given by observing
\begin{equation}\label{EqYRegr}
Y_i=X(i/n)+\eps_i,\quad i=0,1,\ldots,n \text{ with } \eps_i\sim N(0,\eta^2)\text{ i.i.d.}
\end{equation}
for some noise level $\eta>0$. As in classical nonparametrics the statistical analysis becomes much more transparent in equivalent continuous-time observation models, not exposed to discretisation errors. By linear interpolation between the observations $(Y_i)$, \citet{bibinger2014} establish asymptotic equivalence
in Le Cam's sense between the regression model \eqref{EqYRegr} and observing
\begin{equation}\label{EqY}
dY(t)=X(t)dt+\eta n^{-1/2}dW(t), \quad t\in[0,1],
\end{equation}
with a Brownian motion $W$ independent of $X$. In fact, such a result holds for more general observation schemes than \eqref{EqYRegr}, involving possibly non-equidistant and asynchronous observations. For this, $\Sigma(t)$ may also be rank-deficient and no lower bound on the smallest eigenvalue is necessary, as the equivalence proof \citep[Thm. 3.4]{bibinger2014} shows. For the analysis we shall work with model \eqref{EqY}. The transfer to model \eqref{EqYRegr} will be discussed in Section \ref{SecImpl}. The observables in model \eqref{EqY} for test functions $\Phi\in C^1([0,1])$ with $\Phi(0)=\Phi(1)=0$  are given by
\begin{align*}
S_\Phi &:=\int_0^1 \Phi'(t)\,dY(t)=-\int_0^1 \Phi(t)\,dX(t)+\eta n^{-1/2}\int_0^1\Phi'(t)\,dW(t)\\
&\sim N\Big(0,\int_0^1\Sigma(t)\Phi(t)^2dt+\eta^2n^{-1}\norm{\Phi'}_{L^2}^2I_d\Big),
\end{align*}
where integration by parts and It\^o isometry were used. This reveals
an underlying spiked covariance structure, but involving a nonparametric matrix function $\Sigma(t)$. Our test statistics will  combine $S_\Phi$ for local trigonometric functions $\Phi$ at different frequencies, resolving the time variability and profiting from information in higher frequencies.

A by now classical result for this model is that the integrated covariance $\int_0^1\Sigma(t)dt$ can be estimated with rate ${\mathcal O}(n^{-1/4})$ and the spot covariance matrix $\Sigma(t)$ for some fixed $t\in[0,1]$ with rate ${\mathcal O}(n^{-\beta/(4\beta+2)})$ when $t\mapsto\Sigma(t)$ is $\beta$-regular, see \citet{bibinger2014b} for our setting and more recently \citet{FigueroaLopez2024} for a general stochastic volatility setting where $\beta=1/2$. See also Remark \ref{RemRates} and Example \ref{ExRates} below.

We aim at testing the null hypothesis ${\mathcal H}_0: \max_{t\in[0,1]}\rank(\Sigma(t))\le r$. This has been considered by \citet{fissler2017} who construct a test via  matrix perturbation  and prove consistency  under a fixed alternative as $n\to\infty$. When looking at the integer-valued rank, there is no clear concept for local alternatives to quantify the power of rank testing. Yet, writing ${\mathcal H}_0: \max_{t\in[0,1]}\lambda_{r+1}(\Sigma(t))=0$ naturally leads to local alternatives of the form ${\mathcal H}_1(v_n): \max_{t\in[0,1]}\lambda_{r+1}(\Sigma(t))\ge v_n$ and we can determine a detection rate $v_n\downarrow 0$ such that a level-$\alpha$ test under ${\mathcal H}_0$ is asymptotically consistent uniformly over ${\mathcal H}_1(R_nv_n)$ for any sequence $R_n\to\infty$ as $n\to\infty$. The detection rate $v_n$ therefore describes the magnitude of an $(r+1)$-st eigenvalue of the spot covariance matrix, which a test detects by rejecting the null hypothesis of rank at most $r$.

Since we face a nonparametric testing problem, we first concentrate on testing locally ${\mathcal H}_0$ on a sub-interval $I_{t,h}=[t,t+h]$. As a natural test statistic we consider the $(r+1)$st eigenvalue $\lambda_{r+1}(\widehat\Sigma^{t,h})$ of a covariance matrix estimator $\widehat\Sigma^{t,h}$. We construct $\widehat\Sigma^{t,h}$ via blockwise spectral statistics on $I_{t,h}$ as in \citet{bibinger2014b}, but with  adapted
weights, depending on a frequency mixing parameter $M$. We derive explicit and non-asymptotic critical values to ensure level $\alpha$ uniformly over ${\mathcal H}_0$, assuming $\beta$-H\"older regularity and a minimal spectral gap $\underline\lambda_r\ge 0$ for the $r$th eigenvalue of $\Sigma(s)$ on $I_{t,h}$. In contrast to direct observations \citep{RW23} we must accommodate for sums of dependent, non-identically distributed matrix summands, which is achieved by extending classical matrix concentration results \citep{tropp2012} and their treatment via matrix stochastic analysis \citep{bacry2018}. The deviation inequality for the maximal eigenvalue obtained in Theorem \ref{ThmMatrixDevMain} might be of independent interest.
The critical value is primarily determined by two components: a bias term, arising from the time-varying eigenspace and increasing with $h$ (while depending on $\beta$ and $\underline\lambda_r$), and a variance-like term, attributable to observation noise and decaying with both $n$ and $h$. If the bias does not dominate, we obtain asymptotically exact critical values via quantiles of maximal eigenvalues of the Gaussian orthogonal ensemble (GOE).

The power analysis under the alternative ${\mathcal H}_1$ on $I_{t,h}$ shows that the optimal detection rate is governed by the size of the critical value. With an optimal choice of the block length $h$ the local detection rate becomes ${\mathcal O}((\underline\lambda_r^{-1/(\beta+1)}n^{-\beta/(\beta+1)})\wedge n^{-\beta/(\beta+2)})$. In the typical stochastic volatility case $\beta=1/2$ this becomes ${\mathcal O}(n^{-1/3})$ for a positive spectral gap $\underline\lambda_r>0$ and ${\mathcal O}(n^{-1/5})$ without spectral gap. This is much faster than the spot covariance estimation rate ${\mathcal O}(n^{-1/8})$ obtained for $\beta=1/2$.  By multiple testing on all subintervals $I_{kh,h}$ a global test for ${\mathcal H}_0$ on $[0,1]$ is finally obtained whose detection rate is only logarithmically worse than for the local test.
A minimax lower bound shows that the local detection rate $v_n$ is  optimal.

Simulation results confirm the main theoretical findings and also exhibit interesting finite-sample behaviour. An empirical example with intraday bond prices nicely illustrates the impact of the block length $h$ on rank detection. In particular, the rank of the integrated covariance matrix happens to be significantly larger than the maximal rank of the spot covariance matrices. This provides insights into the non-trivial eigenstructure of bond return covariance matrices, a feature recently emphasized by \cite{cg22} in determining the number of factors driving the term structure of interest rates.

In the next Section \ref{SecGeneral} the basic covariance estimators are introduced and a CLT is proved to exhibit the asymptotic error structure in covariance estimation under noise. Section \ref{SecLocalTest} introduces the local tests and analyses their behaviour under the null hypothesis.
The power analysis and derivation of optimal detection rates is presented in Section \ref{SecLocalPower} and the extension to the global test is discussed in Section \ref{SecUMTests}. Section \ref{SecLowerBound} establishes the lower bounds and Section \ref{SecImpl} discusses the implementation, simulation examples and the real data analysis. Longer or more technical proofs are deferred to the Appendix.

\section{Basic results}\label{SecGeneral}

Let us  introduce some  notation. We write $a_n\lesssim b_n$ if $a_n\le Cb_n$ for some constant $C>0$ and all $n$ (or other parameters for the asymptotics). By $a_n\thicksim b_n$ we mean  $a_n\lesssim b_n$ and $b_n\lesssim a_n$. If even $a_n/b_n\to 1$ holds, then we write $a_n\asymp b_n$. With ${\mathcal O}_P(b_n)$ and ${o}_P(b_n)$ we denote random variables $X_n$ such that $b_n^{-1}X_n$, $n\ge 1$, remain bounded in probability and tend to zero in probability, respectively.

We write $v^{\otimes 2}=vv^\top\in\R^{d\times d}$ for vectors $v\in\R^d$. For symmetric matrices $A,B\in\R^{d\times d}$ the partial order $A\le B$ says that $B-A$ is positive semi-definite. Sometimes we also write $A_n\lesssim B_n$ for matrices to say that $A_n\le CB_n$ for all $n$ and some real constant $C>0$. We introduce the set $\R^{d\times d}_{spd}:=\{S\in\R^{d\times d}\,|\, S\text{ symmetric and }S\ge 0\}$. For a symmetric matrix $S\in\R^{d\times d}$ we consider the ordered eigenvalues $\lambda_{max}(S):=\lambda_1(S)\ge\cdots\ge\lambda_d(S)$ (according to their multiplicities), $\trace(S)=\sum_{j=1}^d\lambda_j(S)$ and the remaining traces $\trace_{\ge r}(S)=\sum_{j={r}}^d\lambda_j(S)$, $\trace_{>r}(S)=\sum_{j={r+1}}^d\lambda_j(S)$. For matrices $S,T\in\R^{d\times d}$ the Hilbert-Schmidt or Frobenius scalar product is $\scapro{S}{T}_{HS}=\trace(T^\top S)$ with Hilbert-Schmidt norm $\norm{S}_{HS}=\scapro{S}{S}_{HS}^{1/2}$.  The identity matrix is denoted by $I_d\in\R^{d\times d}$ and the spectral norm is $\norm{S}=\max_{v\in\R^d,\norm{v}=1}\norm{Sv}$. These notions extend naturally to linear operators on Hilbert spaces.

Standard sequence $\ell^p$-norms $\norm{(a_n)}_{\ell^p}=(\sum_n\abs{a_n}^p)^{1/p}$ and $\norm{(a_n)}_{\ell^\infty}=\sup_n\abs{a_n}$ are employed. We use $L^p(I)$ for the standard $L^p$-space of real-valued functions on an interval $I$ with respect to Lebesgue measure. $L^2(I;\R^{d\times d})$ is the Hilbert space of square-integrable matrix functions $S:I\to \R^{d\times d}$ on an interval $I$ with scalar product $\scapro{S_1}{S_2}_{L^2}=\int_I\scapro{S_1(t)}{S_2(t)}_{HS}\,dt$. Finally, for $\beta\in(0,1]$ and $L>0$ we shall consider the matrix-valued H\"older balls with respect to the spectral norm
\[ C^\beta(I,L):=\{\Sigma: I\to\R_{spd}^{d\times d}\,|\, \norm{\Sigma(t)-\Sigma(s)}\le L\abs{t-s}^\beta\text{ for all } t,s\in I\}.
\]

Let us define localised estimators of the spot covariance matrix $\Sigma(t)$, based on the observations of $dY(t)$ from \eqref{EqY}.
We introduce the intervals or blocks
\[I_{t,h}=[t,t+h]\text{ for } t\in[0,1),\,h\in(0,1-t].\]
 On $I_{t,h}$ we form for $j\ge 1$ the local spectral statistics
\begin{align}\label{EqSj}
S_{j}&:=\int_{I_{t,h}} \Phi_{j}'(s)dY(s),\quad \Phi_{j}(s)=\sqrt{2/h}\sin\big(j\pi (s-t)/h\big).
\end{align}
We omit the dependence of $\Phi_j$ on $t$ and $h$ for the sake of brevity. Note that $\{\Phi_j\,|\,j\ge 1\}$ and $\{\Phi_j'/\norm{\Phi_j'}_{L^2(I_{t,h})}\,|\,j\ge 1\}$ form orthonormal systems in $L^2(I_{t,h})$. Since $\Phi_j(t)=\Phi_j(t+h)=0$ holds, we obtain by partial integration
\begin{equation}\label{EqSj2}
 S_{j}=-\int_{I_{t,h}} \Phi_j(s)\Sigma^{1/2}(s)\,dB(s)+\int_{I_{t,h}} \Phi_j'(s)\eta n^{-1/2}\,dW(s).
 \end{equation}

By independence of $B$, $W$, using It\^o isometry and $2\sin^2(\alpha)=1-\cos(2\alpha)$, this gives $S_{j}\sim N(0,C_{j})$ with
\begin{align}
C_{j}&=\int_{I_{t,h}}\Sigma(s)\Phi_j(s)^2ds+\norm{\Phi_j'}_{L^2(I_{t,h})}^2\eta^2n^{-1}I_d
=\Sigma_{j}+j^2\eps_{n,h}^2I_d,\text{ where}\label{EqCj}\\
 \Sigma_{j}&:=\frac1h\int_{I_{t,h}}\Sigma(s)\Big(1
-\cos\big(2\pi j(s-t)/h\big)\Big)\,ds,\\
\eps_{n,h} &:=\pi h^{-1}\eta n^{-1/2}.\label{Eqeps}
\end{align}
We call $\eps_{n,h}$ the local noise level or inverse signal-to-noise ratio in \eqref{EqCj}. We can think of $nh$ observations in $I_{t,h}$ for the regression model \eqref{EqYRegr} with errors of size $\eta$, while the signal $\Sigma(s)^{1/2}$ on $I_{t,h}$ is multiplied by $h^{1/2}$, the size of the Brownian increment on $I_{t,h}$, so that the signal-to-noise ratio on $I_{t,h}$ scales like $h^{1/2}/(\eta^2/nh)^{1/2}\thicksim\eps_{n,h}^{-1}$. Even though we shall profit from averaging over different frequencies $j$, a minimal block size $h\gtrsim n^{-1/2}$ will become necessary, that is $\eps_{n,h}\lesssim 1$ in \eqref{Eqeps}.

We obtain for $j'\not=j$ by orthogonality of $\Phi_j'$ and $\Phi_{j'}'$
\begin{equation}\label{EqSjSj'} \E[S_{j}S_{j'}^\top]=\int_{I_{t,h}} \Sigma(s)\Phi_j(s)\Phi_{j'}(s)\,ds.
\end{equation}
This implies independence of $S_{j}$, $S_{j'}$ when $\Sigma(s)$ is constant on $I_{t,h}$. The key idea of the localized spectral statistics is that, for regular functions $\Sigma(s)$ and small block size $h$, they exhibit only weak correlation across different frequencies $j$.

For $M\ge 1$ we consider weights, which lead to optimised block-wise estimation results in the critical case $\Sigma(t)=M^2\eps_{n,h}^2I_d$, compare the oracle weights in \citet{bibinger2014} and also Example \ref{ExRates} below:
\begin{equation}\label{Defwj}
 w_j:=c_w M^{-1}(1+j^2/M^2)^{-2},\quad j\ge 1,
\end{equation}
with (asymptotics as $M\to \infty$)
\[ c_w:=\Big(\sum_{j\ge 1}M^{-1}(1+j^2/M^2)^{-2}\Big)^{-1}\asymp \Big(\int_0^\infty(1+x^2)^{-2}dx\Big)^{-1}=\frac4{\pi}
\]
and $c_w\le 4/\pi$ for all $M$ (use the integral bound). The dependence of $w_j$, $c_w$ on $M$ is omitted in the notation.
We obtain an estimator of $\Sigma$ on the interval $I_{t,h}$ via
\begin{align}\label{EqhatSigma}
 \widehat\Sigma^{t,h} &:=\widehat C^{t,h}-B_w\eps_{n,h}^2I_d,\quad\text{where }
 \widehat C^{t,h}:=\sum_{j\ge 1} w_j S_{j}^{\otimes 2},\,B_w:=\sum_{j\ge 1} w_jj^2.
\end{align}
The bias correction with $B_w$ is chosen such that in terms of $C_j$, $\Sigma_j$ from \eqref{EqCj}
\begin{align}
C^{t,h} &:=\E[\widehat C^{t,h}]=\sum_{j\ge 1}w_j C_{j},\label{EqCth}\\
\Sigma^{t,h} &:=\E[\widehat\Sigma^{t,h}]=\sum_{j\ge 1}w_j \Sigma_{j}
=\int_0^1 \Sigma(t+hs)w(s)\,ds,\label{EqSigmath}\\
\text{for } w(s)&:=1-\sum_{j\ge 1} w_j\cos\big(2\pi js\big).\nonumber
\end{align}
Observe that the weight function satisfies $w(s)\ge 0$ and integrates to one so that $\Sigma^{t,h}$ is a weighted time average of $\Sigma$ on $I_{t,h}$.
The weights satisfy $w_j\thicksim M^{-1}$ for $j\lesssim M$ and decay rapidly for $j\gg M$ so that in practice $\widehat C^{t,h}$ is obtained as a convex combination over ${\mathcal O}(M)$ empirical covariance matrices $S_{j}S_{j}^\top$. The larger $M$, the more we average  over frequencies, and we call $M$ the { mixing parameter}.

We start our analysis by showing asymptotic normality of $\widehat\Sigma^{t,h}$ and exhibiting the influence of the sample size $n$, the block length $h$, the mixing parameter $M$ and the covariance matrix itself on the convergence rate.

As a preparation let us recall that the Gaussian orthogonal ensemble $GOE(d)$ on the symmetric $d\times d$-matrices is given by the law of $Z=(Z_{i,j})_{i,j=1,\ldots,d}$ with independent $Z_{i,j}$ for $j\ge i$, $Z_{j,i}=Z_{i,j}$ and $Z_{i,i}\sim N(0,2)$, $Z_{i,j}\sim N(0,1)$ for $j>i$. We use the $\vek$-operation that stacks a $d\times d$-matrix into a vector in $\R^{d^2}$ and the Kronecker product $\otimes$ between matrices, see e.g. Section 2.2 in \cite{bibinger2014}. Then in view of \eqref{EqSjSj'}
the covariance structure of the covariance matrix estimator in terms of  ${\mathcal Z}_d=\Cov(\vek(Z^{\otimes 2}))\in\R^{d^2\times d^2}$ is given by
\begin{align*}
\Cov(\vek(\widehat\Sigma^{t,h}))&=\Cov\Big(\sum_{j\ge 1}w_j\vek(S_j^{\otimes 2})\Big)={\mathcal C}(\widehat\Sigma^{t,h}){\mathcal Z}_d\\
\text{with }{\mathcal C}(\widehat\Sigma^{t,h})&:=\sum_{j\ge 1} w_j^{2}  C _{j}^{\otimes 2}+\sum_{j\not=j'} w_jw_{j'}\Big(\int_{I_{t,h}}\Sigma(s)\Phi_j(s)\Phi_{j'}(s)\,ds\Big)^{\otimes 2}.
\end{align*}

\begin{theorem}\label{ThmCLT}
Consider the  asymptotics $M\to\infty$, while $t,h,n$ may vary arbitrarily with $M$. For each $v\in\R^d$ with $\norm{v}=1$ assume
\begin{equation}\label{EqAssCLT}
\exists C_v>0\,\forall M\ge 1:\; \norm{\scapro{\Sigma(\cdot)v}{v}}_{L^\infty(I_{t,h})}\le C_v \Big( M^2\eps_{n,h}^2+h^{-1}\int_{I_{t,h}}\scapro{\Sigma(s)v}{v}\,ds\Big).
\end{equation}
Then the following central limit theorem  holds:
\[ {\mathcal C}(\widehat\Sigma^{t,h})^{-1/2}\vek\big(\widehat\Sigma^{t,h}-\Sigma^{t,h}\big)\xrightarrow{d} N(0,{\mathcal Z}_d)=\vek(GOE(d)).
\]
Moreover, we have ${\mathcal C}(\widehat\Sigma^{t,h})\lesssim M^{-1}((\norm{\Sigma_{\ell,\ell'}}_{L^\infty(I_{t,h})})_{\ell,\ell'=1,\ldots,d}+\eps_{n,h}^2M^2I_d)^{\otimes 2}$.
\end{theorem}

\begin{proof} See Appendix \ref{AppSecProof2}. \end{proof}

The assumption in \eqref{EqAssCLT} is satisfied if the local noise level dominates or if the average and maximum of $\scapro{\Sigma(\cdot)v}{v}$ on $I_{t,h}$ are comparable, which only excludes very spiky behaviour.
The proof of the CLT is considerably simplified by the fact that $\widehat\Sigma^{t,h}$  belongs to a second Wiener chaos.

The result implies that $\widehat\Sigma^{t,h}$  estimates for fixed $t$ and $h\downarrow 0$ the spot covariance matrix $\Sigma(t)$ consistently, provided $\eps_{n,h}^4M^3\to 0$ and $\Sigma$ is continuous at $t$. For fixed $M$ the estimator $\widehat\Sigma^{t,h}$ will not be consistent for  $\Sigma^{t,h}$ even in the extreme case  $\eps_{n,h}=0$ of no additive noise because it is a fixed convex combination of covariance estimators whose error has the size of the covariance itself.

\begin{remark}\label{RemRates}
 In related work on estimation of the integrated covariance \citep{ait2010,barndorff2011,christensen2013,bibinger2014}, the methods are tuned for covariance matrices  with smallest and largest eigenvalues of order one. In our setting, this leads to the choice $M\thicksim \eps_{n,h}^{-1}$ resulting in a blockwise convergence rate $\eps_{n,h}^{1/2}\thicksim h^{-1/2}n^{-1/4}$ for $\widehat\Sigma^{t,h}$ and at best an $n^{-1/4}$-rate for the integrated covariance $\int_0^1 \Sigma(s)ds$ by averaging over $\floor{h^{-1}}$ disjoint blocks on $[0,1]$. For testing zero eigenvalues of the covariance, however, we shall choose $M$ much smaller (or even keep it fixed) so that we attain almost the convergence rate $\eps_{n,h}^2\thicksim h^{-2}n^{-1}$ for the small empirical eigenvalues in $\widehat\Sigma^{t,h}$. The loss in precision for the larger empirical eigenvalues does not spoil the test statistics and we shall establish a much faster detection rate than the usual $n^{-1/4}$-estimation rate would suggest.
\end{remark}

\section{\bf Local tests  under the null hypothesis}\label{SecLocalTest}

We aim to test whether the covariance function $\Sigma(s)$ has at most rank $r$ on a given interval $I_{t,h}$, where $r\in\{0,\ldots,d-1\}$ is fixed throughout. The natural test statistics is the $(r+1)$st eigenvalue of $\widehat C^{t,h}$ (or equivalently $\widehat\Sigma^{t,h}$). Our power analysis will establish a sort of bias-variance dilemma: for a small length $h$ of the observation interval the stochastic error is large while averaging with large $h$ can induce large $\lambda_{r+1}(\Sigma^{t,h})$ even though $\lambda_{r+1}(\Sigma(s))=0$ pointwise on the block. Later, we will determine an asymptotically optimal $h$ that maximises power under local alternatives. By a lower bound argument we demonstrate that the bias-variance trade-off is not an artifact of our specific spectral covariance matrix estimator, but an intrinsic challenge in rank testing of time-varying covariance matrices.

First let us define precisely the composite null hypothesis we want to test.

\begin{definition}
Let  $\beta\in(0,1]$, $\underline\lambda_r\ge 0$, $L>0$. The null hypothesis of at most rank $r$ on $I_{t,h}$ is given by
\begin{align*}
{\mathcal H}_0&:={\mathcal H}_0(I_{t,h},\beta,L,\underline\lambda_{r})\\
&:=\big\{\Sigma\in C^\beta(I_{t,h},L)\,|\,\rank(\Sigma(s))\le r,\,\lambda_r(\Sigma(s))\ge\underline\lambda_{r}\text{ for all $s\in I_{t,h}$}\big\}.
\end{align*}
\end{definition}

Besides the H\"older regularity constraint expressed by $\beta$ and $L$, we also impose a spectral gap condition  in ${\mathcal H}_0$ on the $r$th eigenvalue of the covariance matrix in terms of $\underline\lambda_r$. In the case $\underline\lambda_r=0$ there is no constraint, but we shall be able to gain in the testing power by a positive spectral gap.

In order to construct a test based on the estimator $\widehat C^{t,h}$, let us express its expectation $C^{t,h}$ from \eqref{EqCth} in an orthonormal basis of eigenvectors $v_1,\ldots,v_d$ with eigenvalues $\lambda_1( C^{t,h})\ge\cdots\ge\lambda_d( C^{t,h})\ge 0$. We introduce the space $V_{>r}:=\spann(v_{r+1},\ldots,v_d)$ and the orthogonal projection $P_{>r}$ onto $V_{>r}$ as well as the $(d-r)\times (d-r)$-lower right minors
\[ S_{>r}:=P_{>r}S|_{V_{>r}}\text{ for matrices } S\in\{\Sigma(s),\widehat \Sigma^{t,h},\widehat C^{t,h}, \Sigma^{t,h}, C^{t,h}\}.\]
Our analysis uses repeatedly that by the Cauchy interlacing law (e.g., \citet{johnstone2001} or \citet{tao2012}), the  estimators  $\widehat \Sigma^{t,h}$, $\widehat C^{t,h}$ satisfy
\begin{equation}\label{EqCIL}
\lambda_{r+1}(\widehat \Sigma^{t,h})\le \lambda_{max}(\widehat \Sigma^{t,h}_{>r}),\quad \lambda_{r+1}(\widehat C^{t,h})\le \lambda_{max}(\widehat C^{t,h}_{>r}).
\end{equation}

For $\Sigma(s)$ constant on $I_{t,h}$ we would have $\Sigma^{t,h}_{>r}=0$ leading to a  spiked covariance model for $\widehat C^{t,h}$ \citep{johnstone2001}. In general, though, $\lambda_{r+1}(\Sigma^{t,h})>0$ holds under ${\mathcal H}_0$, but its  magnitude
is bounded by
\begin{equation}\label{EqBiasRW}
\norm{\Sigma_{>r}(\cdot)}_{L^\infty(I_{t,h})}\le (2Lh^\beta) \wedge ((r+4)\tfrac{L^2h^{2\beta}}{\underline\lambda_r})=:\BIAS_0^\Sigma(h),
\end{equation}
which follows from Prop. S.4(c) in \citet{RW23} (considering $p=\infty$ in that bound, $I'=[s,s+\delta]$ for $s\in I_{t,h}$ and letting $\delta\downarrow 0$).

A second key ingredient is an upper deviation inequality for the maximal eigenvalue of the matrix $\sum_{j\ge 1}\Gamma_j^{\otimes 2}$ with Gaussian random vectors $\Gamma_j$ involving different covariance matrices and arbitrary correlation. This is essential for bounding $\lambda_{max}(\widehat C^{t,h}_{>r})$  because the spectral statistics $(S_j)$ have an intricate variance-covariance structure, which cannot be overcome by a locally constant approximation of $\Sigma$ under the H\"older assumption. The proof is accomplished via tools from the Bernstein inequality for continuous matrix martingales and via Gaussian concentration.

\begin{theorem}\label{ThmMatrixDevMain}
For $G_j\in L^2([0,1];\R^{d\times d})$ consider $\Gamma_j=\int_0^1G_j(u)\,dB(u)$, $j\ge 1$, with a $d$-dimensional Brownian motion $B$ and $\sum_{j\ge 1}\norm{G_j}_{L^2}^2<\infty$. Then for $\alpha\in(0,1)$ we have with probability at least $1-\alpha$
\begin{align*}
\lambda_{max}\Big(\sum_{j\ge 1}\Gamma_j^{\otimes 2}\Big)\le \inf_{\delta>0}\Big(&(1+\delta)\Big(\lambda_{max}\Big(\sum_{j\ge 1} \E[\Gamma_j^{\otimes 2}] \Big)+ \log(\sqrt{7e/2}d)\max(v_1,v_2)\Big)\\
& +(1+\delta^{-1})2\log(\alpha^{-1}) \sigma^2\Big),
\end{align*}
where
\begin{align*}
v_1^2 &=\bnorm{\E\Big[\Big(\sum_{j\ge 1}\Gamma_j^{\otimes 2}\Big)^2\Big]-\Big(\E\Big[\sum_{j\ge 1}\Gamma_j^{\otimes 2}\Big]\Big)^2},\\
v_2^2 &=\sup_{\norm{F}_{L^2([0,1];\R^{d\times d})}\le 1}\bnorm{\sum_{j\ge 1}\int_\cdot^1\big(\scapro{G_j(s)}{F(s)}_{HS}I_d+G_j(s)F(s)^\top\big)\, ds\,G_j}_{L^2([0,1];\R^{d\times d})}^2,\\
\sigma^2&= \sup_{\norm{f}_{L^2([0,1];\R^d)}\le 1}\bnorm{\sum_{j\ge 1}\Big(\int_0^1 G_jf\Big)^{\otimes 2}}.
\end{align*}
\end{theorem}

\begin{proof}
This reformulates Theorem \ref{ThmMatrixDev} below.
\end{proof}

\begin{example}
The purpose of Theorem \ref{ThmMatrixDevMain} is to cover anisotropic, dependent random vectors, but
let us check the result in the simple, yet later useful case of independent $\Gamma_j\sim N(0,s_j I_d)$, $s_j\ge 0$, with $\sum_js_j<\infty$. Then
\[ \lambda_{max}\Big(\sum_{j\ge 1}\E[\Gamma_j^{\otimes 2}]\Big)=\sum_{j\ge 1}s_j,\quad v_1^2=\sum_{j\ge 1} (d+1)s_j^2
\]
follows directly. The expressions for $v_2$ and $\sigma$ simplify considerably if we choose $G_j(t)=s_j^{1/2}f_j(t)I_d$  for some orthonormal system $(f_j)$ in $L^2([0,1])$ with disjoint support (e.g., $f_j=(b_j-a_j)^{-1/2}{\bf 1}_{[a_j,b_j]}$ for pairwise disjoint intervals $[a_j,b_j]$). By the Bessel inequality we then have $\sum_j(\int_0^1f_jf)^2\le \norm{f}_{L^2}^2$ so that
$\sigma^2\le \max_{j\ge 1}s_j$
follows. The isotropy of $G_j$ in the coordinates leads to a maximiser of the form $F(t)=d^{-1/2}g(t)I_d$ with $\norm{g}_{L^2([0,1];\R)}=1$ when bounding via Cauchy-Schwarz inequality
\begin{align*}
v_2^2 &= \sup_{\norm{F}_{L^2}\le 1}\sum_{j\ge 1}s_j^2\int_0^1 \bnorm{\int_t^1\big(\scapro{I_d}{F(s)}_{HS}I_d+F(s)^\top\big)f_j(s)\, ds}_{HS}^2f_j(t)^2 dt\\
&\le \sup_{\norm{g}_{L^2}\le 1}\sum_{j\ge 1}s_j^2\int_0^1 (d+1)^2d\norm{d^{-1/2}g}_{L^2(\supp(f_j))}^2\norm{f_j}_{L^2}^2f_j(t)^2 dt\\
&=(d+1)^2\sup_{j\ge 1}s_j^2.
\end{align*}
Hence, we obtain with probability at least $1-\alpha$
\begin{align}\label{EqDevIneqEx}
&\lambda_{max}\Big(\sum_{j\ge 1}\Gamma_j^{\otimes 2}\Big)
 \le \inf_{\delta>0}\Big((1+\delta)\norm{(s_j)}_{\ell^1}+ (1+\delta)\log(\sqrt{7e/2}d)\\
 &\quad\times\max\big(\sqrt{d+1}\norm{(s_j)}_{\ell^2},(d+1)\norm{(s_j)}_{\ell^\infty}\big) +2(1+\delta^{-1})\log(\alpha^{-1})\norm{(s_j)}_{\ell^\infty}\Big).\nonumber
\end{align}
Neglecting logarithmic factors this bound is of order $\norm{(s_j)}_{\ell^1}\vee d\norm{(s_j)}_{\ell^\infty}$, which is easily seen to be tight. Compare also the  analogous bound derived by the matrix Bernstein inequality under independence \cite[Thm. A.2]{RW23}.
\end{example}

We are now in a position to define and analyse a  non-asymptotic test, which rejects ${\mathcal H}_0$ if $\lambda_{r+1}(\widehat C^{t,h})$ is large. For this recall $\BIAS_0^\Sigma(h)$ from \eqref{EqBiasRW}.

\begin{theorem}\label{Propblockwisetest0}
Fix $\alpha\in(0,1)$, $M\ge r$ and define the test $\phi_{\alpha,0}={\bf 1}(\lambda_{r+1}(\widehat C^{t,h})>\kappa_{\alpha,0})$  with critical value
\begin{align}
\kappa_{\alpha,0}&:=
  \big(4\BIAS_0^\Sigma(h)+\eps_{n,h}^2\tfrac{M^2}2\big) \Big(1+
 \log(4(d-r)){\mathcal W}\big(\tfrac{6d}{M}\big)
 +\tfrac{21\log(2\alpha^{-1})}M \Big),
 \label{Eqkappa0}
\end{align}
where ${\mathcal W}(x):=x\vee\sqrt{x}$, $x\ge 0$.
Then $\phi_{\alpha,0}$  is a non-asymptotic uniform level-$\alpha$ test for ${\mathcal H}_0$:
\[ \sup_{\Sigma\in{\mathcal H}_0(I_{t,h},\beta,L,\underline\lambda_r)}\PP_{\Sigma}(\phi_{\alpha,0}=1)\le \alpha.\]
\end{theorem}

\begin{proof}
See Corollary \ref{CorChat} below.
\end{proof}

\begin{remark}
The  critical value $\kappa_{\alpha,0}$ reflects the usual multiplicative structure in the stochastic error of covariance estimators.
Larger choices of $M$ reduce the stochastic error by averaging better over the frequencies. On the other hand, the  noise part induced by $\eps_{n,h}$
grows in $M$.
The numerical constants in \eqref{Eqkappa0} are explicit, but not  optimised.
\end{remark}

In the next result we  provide asymptotically tight critical values in the case where the bias does not dominate.

\begin{theorem}\label{Thmblockwisetest}\mbox{}
Consider the asymptotics $M\to \infty$ and $\BIAS_0^\Sigma(h)\lesssim \eps_{n,h}^2$, where $n,h,\underline\lambda_r$ can vary arbitrarily otherwise. Then
 \begin{align}
 \phi_{\alpha,1}&:={\bf 1}(\lambda_{r+1}(\widehat C^{t,h})>\kappa_{\alpha,1})\text{ with}\\
\kappa_{\alpha,1}&:=\big(B_w+(2\pi)^{-1/2}M^{3/2}q_{1-\alpha;\lambda_{max}(GOE(d-r))}\big) \eps_{n,h}^2
\end{align}
and the $(1-\alpha)$-quantile $q_{1-\alpha;\lambda_{max}(GOE(d-r))}$ of $\lambda_{max}(GOE(d-r))$
is an asymptotically uniform level-$\alpha$ test:
\[ \limsup_{M\to\infty}\sup_{\Sigma\in{\mathcal H}_0(I_{t,h},\beta,L,\underline\lambda_r)} \PP_\Sigma(\phi_{\alpha,1}=1)\le\alpha.
\]
This remains true for $\phi_{\alpha,2}:={\bf 1}(\lambda_{r+1}(\widehat C^{t,h})>\kappa_{\alpha,2})$ with the simulation-based critical value $\kappa_{\alpha,2}:=q_{1-\alpha,\Lambda_M}\eps_{n,h}^2$ where $q_{1-\alpha,\Lambda_M}$ is the $(1-\alpha)$-quantile of
\begin{equation}\label{eqLambdaM}
 \Lambda_M:=\lambda_{max}\Big(\sum_{j\ge 1}w_jj^2\zeta_j^{\otimes 2}\Big)\text{ with independent } \zeta_j\sim N(0,I_{d-r}).
 \end{equation}
\end{theorem}

\begin{proof}
Representing $S_{j}=\Sigma_{j}^{1/2}Z_j+\eps_{n,h}^2j^2Z_j'$ with  $Z_j,Z_j'\sim N(0,I_d)$, we can split $\widehat C^{t,h}=\widehat C^{t,h}_1+\widehat C^{t,h}_2+\widehat C^{t,h}_3$ with
\begin{align}
\widehat C^{t,h}_1 &= \sum_{j\ge 1}w_j(\Sigma_{j}^{1/2}Z_j)^{\otimes 2},\qquad
\widehat C^{t,h}_2 = \eps_{n,h}^2\sum_{j\ge 1}w_jj^2(Z_j')^{\otimes 2},\label{EqDecompChat}\\
\widehat C^{t,h}_3 &= \eps_{n,h}\sum_{j\ge 1}w_jj\Big(\Sigma_{j}^{1/2}Z_j(Z_j')^\top + Z_j'Z_j^\top\Sigma_{j}^{1/2}\Big).\nonumber
\end{align}
From \eqref{EqCIL} we infer
\[\lambda_{r+1}(\widehat C^{t,h})\le \lambda_{max}(\widehat C^{t,h}_{1,>r})+ \lambda_{max}(\widehat C^{t,h}_{2,>r})+\lambda_{max}(\widehat C^{t,h}_{3,>r}).
\]
By Corollary \ref{CorChat} below, inserting $\eta=0$ (no noise) and taking any sequence $\alpha_M\to 0$ with $\log(\alpha_M^{-1})M^{-1}\to 0$,  the first term $\widehat C^{t,h}_{1,>r}$ satisfies
\begin{align*}
\lambda_{max}\big(\widehat C^{t,h}_{1,>r}\big)&={\mathcal O}_P(\BIAS_0^\Sigma(h))
={\mathcal O}_P( \eps_{n,h}^2),
\end{align*}
uniformly over $\Sigma\in{\mathcal H}_0$ in view of $\BIAS_0^\Sigma(h)\lesssim \eps_{n,h}^2$.
The CLT from Theorem \ref{ThmCLT}, applied in the pure noise case ($\Sigma\equiv 0$) and dimension $d-r$, yields  that
\[ \rho_2^{-1/2}\big(\widehat C^{t,h}_{2,>r}-\E[\widehat C^{t,h}_{2,>r}]\big)\xrightarrow{d} GOE(d-r)\text{ as } M\to\infty\]
with $\rho_2=\norm{j^2w_j}_{\ell^2}^2 \eps_{n,h}^4\asymp \tfrac1{2\pi}M^3 \eps_{n,h}^4$ and $\E[\widehat C^{t,h}_{2,>r}]=B_w\eps_{n,h}^2I_{d-r}$, where all other quantities $t,h,n,\underline\lambda_r$ may vary arbitrarily.
Since the eigenvalue map $S\mapsto\lambda_{max}(S)$ is continuous on the set of symmetric matrices $S$, the continuous mapping theorem implies
\[ \sqrt{2\pi}M^{-3/2}\big(\eps_{n,h}^{-2}\lambda_{max}\big(\widehat C^{t,h}_{2,>r}\big)-B_w\big) \xrightarrow{d} \lambda_{max}\big(GOE(d-r)\big).
\]
For the cross term we obtain by the Cauchy-Schwarz inequality
\[ \norm{\widehat C_3^{kh}}\le 2\eps_{n,h}\sum_{j\ge 1}w_jj\norm{\Sigma_{jk}^{1/2}Z_j}\norm{Z_j'} \le 2\trace(\widehat C_1^{kh})^{1/2}\trace(\widehat C_2^{kh})^{1/2},
\]
which by the above results is ${\mathcal O}_P(\eps_{n,h}){\mathcal O}_P(\eps_{n,h}M)=o_P(M^{3/2}\eps_{n,h}^2)$.

Consequently, Slutsky's Lemma applies and yields for all $q\in\R$
\[ \limsup_{M\to\infty}\PP\Big(\sqrt{2\pi}M^{-3/2} \big(\eps_{n,h}^{-2}\lambda_{r+1}(\widehat C^{t,h})-B_w\big)\ge q\Big)\le \PP\Big(\lambda_{max}(GOE(d-r))\ge q\Big)\]
uniformly over $\Sigma\in {\mathcal H}_0$, which by using the quantile $q=q_{1-\alpha;\lambda_{max}(GOE(d-r))}$ gives the result.

Concerning $\kappa_{\alpha,2}$, we just remark that $\eps_{n,h}^2\Lambda_M$ has the law of $\lambda_{max}(\widehat C_{2,>r}^{kh})$ and we have seen that $\widehat C_{1,>r}^{kh}$ and $\widehat C_{3,>r}^{kh}$ become negligible in norm as $M\to\infty$.
\end{proof}

Note that $\phi_{\alpha,1}$ can be expressed equivalently  by the spot covariance estimator $\widehat\Sigma^{t,h}$:
 \begin{align*}
 \phi_{\alpha,1}&:={\bf 1}\Big(\lambda_{r+1}(\widehat \Sigma^{t,h})>(2\pi)^{-1/2}M^{3/2}q_{1-\alpha;\lambda_{max}(GOE(d-r))} \eps_{n,h}^2\Big).
\end{align*}
The same applies to $\phi_{\alpha,2}$.
The inequality $\lambda_{r+1}(\widehat\Sigma^{t,h})\le \lambda_{max}(\widehat\Sigma^{t,h}_{>r})$ leads to a conservative test, but it becomes tight asymptotically when $\lambda_r(\Sigma^{t,h})$ tends to infinity because the eigenvalue repulsion becomes negligible. This observation shows that the test $\phi_{\alpha,1}$  exhausts the level $\alpha$ on the composite null hypothesis, i.e. the maximal type I error probability over ${\mathcal H}_0$ even tends to $\alpha$ exactly. After Proposition 1.2 in \citet{johnstone2001} the empirical evidence is reported that for large, but finite $\lambda_r(\Sigma^{t,h})$ the law of $\lambda_{r+1}(\widehat\Sigma^{t,h})$ is approximately shifted by a constant from the law of $\lambda_{max}(\widehat\Sigma^{t,h}_{>r})$. A corresponding data-driven adjustment of critical values can be of interest, but so far it lacks theoretical underpinning.

The critical values in the test depend on the remaining dimension $d-r$. For increasing $d$ we know that $(d-r)^{1/6}(\lambda_{max}(GOE(d-r))-\sqrt{2(d-r)})$ approaches the Tracy-Widom law \citep{johnstone2012}. The critical value $\kappa_{\alpha,1}$  will therefore grow like $\sqrt{d-r}$ in the dimension. For a truly high-dimensional setting further assumptions could be imposed, leading to constraint eigenvalue estimation. Compare \citet{ait2017} for a sparse PCA approach based on the integrated covariance matrix, which, however, is not easily transferred to spot covariance matrices  because of the bias incurred by local smoothing.

\section{\bf Power against local alternatives and detection rate}\label{SecLocalPower}

Let us define precisely the local alternatives against which we shall ensure asymptotic power.

\begin{definition}
For $\underline\lambda_{r+1}> 0$ the  alternative hypotheses on $I_{t,h}$ are given by
\begin{align*}
{\mathcal H}_1&:={\mathcal H}_1(I_{t,h},\underline\lambda_{r+1})
:= \Big\{\Sigma\in L^1([0,1];\R_{spd}^{d\times d})\,|\, \lambda_{r+1}(\Sigma(s))\ge \underline{\lambda}_{r+1} \text{, $s\in I_{t,h}$}\Big\}.
\end{align*}
\end{definition}

Note that ${\mathcal H}_1$ does neither involve a H\"older regularity constraint nor a spectral gap condition on $\Sigma$ ($\Sigma\in L^1$ is assumed to define the stochastic integral). The main interest is to establish asymptotic power of our  tests $\phi_{\alpha,i}$ for rank $r$ under  alternatives ${\mathcal H}_1(I_{t,h},\underline\lambda_{r+1})$ which localise around ${\mathcal H}_0$ in the sense that $\underline{\lambda}_{r+1}\downarrow 0$ is considered as $\eps_{n,h}\to 0$. The faster $\underline{\lambda}_{r+1}\downarrow 0$, the better the tests must discriminate zero from small eigenvalues.  First, we provide a general  result, then we optimise the block length $h$.

\begin{theorem}\label{ThmH1blockwise}
Fix $\alpha\in(0,1)$ and $i\in\{0,1,2\}$. The tests $\phi_{\alpha,i}$  are consistent under local alternatives ${\mathcal H}_1(I_{t,h},\underline{\lambda}_{r+1})$ whenever
$M\kappa_{\alpha,i}/\underline{\lambda}_{r+1}\to 0$.
More precisely, we have uniform local consistency
\begin{equation}\label{EqH1blockwise}
\lim_{R\to\infty}\inf_{\Sigma\in {\mathcal H}_1(I_{t,h},RM\kappa_{\alpha,i})} \PP_{\Sigma}(\phi_{\alpha,i}=1)=1,
\end{equation}
where  the parameters $n$, $h$ and $M\ge r+1$ can vary arbitrarily as $R\to\infty$.
\end{theorem}

\begin{proof}
This is Corollary \ref{CorLambdaLower} below.
\end{proof}

The proof refines the analysis in Section S.3 of \citet{RW23}, taking care of the dependence of the spectral statistics $S_j$.
As explained in Remark \ref{RemH1beta} below, we can get rid of the factor $M$, that is consider ${\mathcal H}_1(I_{t,h},R\kappa_{\alpha,i})$ in \eqref{EqH1blockwise}, in some cases, but do not achieve this in general.

We see that the tests are the more powerful, the smaller the critical values $\kappa_{\alpha,i}$ are. For the following asymptotic result, we focus on the nonasymptotic test $\phi_{\alpha,0}$, keep the mixing parameter $M$ fixed and choose the blocksize $h$ to balance the two summands in $\kappa_{\alpha,0}$, thus balancing the stochastic and bias-type errors. For the asymptotic tests $\phi_{\alpha,1}$, $\phi_{\alpha,2}$ the results are only marginally worse when the convergence $M=M_n\to\infty$ is slow for $n\to\infty$.

\begin{corollary}\label{CorBlockTest}
Let $M\ge r+1$  and $\alpha\in(0,1)$  be fixed and choose $h\to 0$ as $n\to\infty$ such that under ${\mathcal H}_0(I_{t,h},\beta,L,\underline\lambda_{r})$
\[\BIAS_0^\Sigma(h)\thicksim \eps_{n,h}^2,\quad\text{i.e. }\quad h\thicksim \big(\underline\lambda_{r}^{-1}L^2n\big)^{-1/(2\beta+2)}\vee  (Ln)^{-1/(\beta+2)}.\]
Then the  level-$\alpha$ tests $\phi_{\alpha,0}$ of the null hypothesis
${\mathcal H}_0(I_{t,h},\beta,L,\underline\lambda_{r})$ satisfy for any sequence $R_n\to\infty$
\[ \lim_{n\to\infty}\inf_{\Sigma\in {\mathcal H}_1(I_{t,h},R_nv_n)} \PP_{\Sigma}(\phi_{\alpha,0}=1)=1
\]
with {\it detection rate}
\begin{equation}\label{Eqvn}
v_n=\big((\underline\lambda_{r}^{-1}L^2)^{1/(\beta+1)}n^{-\beta/(\beta+1)}\big)\wedge  \big(L^{2/(\beta+2)} n^{-\beta/(\beta+2)}\big).
\end{equation}
As $n\to\infty$ the  parameters $\underline\lambda_r$ and $L$ may vary arbitrarily with $n$.
\end{corollary}

\begin{proof}
By the choice of $h$ and the fixed size of $M$ we have
$\kappa_{\alpha,0}\thicksim \BIAS_0^\Sigma(h)+\eps_{n,h}^2\thicksim v_n$ and
 Theorem \ref{ThmH1blockwise} gives the result.
\end{proof}

\begin{example}\label{ExRates}
As a test case we consider $\beta=1/2$ and $L\thicksim 1$  which mimics stochastic volatility models \cite[Section 4]{RW23}. Choosing $h\thicksim (n/\underline\lambda_r)^{-1/3}\vee n^{-2/5}$, the detection rate becomes $v_n\thicksim(\underline\lambda_r^{-2/3}n^{-1/3})\wedge n^{-1/5}$.
Remarkably, for a fixed spectral gap $\underline\lambda_r>0$  our block size becomes $h\thicksim n^{-1/3}$ which is exactly the window size obtained by \citet{fissler2017} with $k_n=nh\thicksim n^{2/3}$ observations in a block. There it was conjectured to be a suboptimal choice dictated by the method employed. For our local test, it is seen to balance the bias of order $h$ for $\beta=1/2$ with the blockwise noise level $\eps_{n,h}^2\thicksim h^{-2}n$. 

The classical estimation rate under noise is $n^{-\beta/(4\beta+2)}$ for spot volatility  \citep{bibinger2014b}. The latter rate is $n^{-1/8}$ for $\beta=1/2$ and our detection rate $(\underline\lambda_r^{-2/3}n^{-1/3})\wedge n^{-1/5}$ seems at first glance surprisingly fast. Note, however, that the estimation rates are stated for fixed covariances $\Sigma$. For $\norm{\Sigma}\to 0$ the spot estimation rate can be found to be $(\norm{\Sigma}^3n^{-1})^{\beta/(4\beta+2)}$. When $\norm{\Sigma}\le n^{-\beta/(\beta+2)}$, we achieve the optimal spot volatility estimation rate $n^{-\beta/(\beta+2)}$, which matches our (slower) rate in the absence of a spectral gap. In contrast, the faster rate observed in the presence of a spectral gap has no direct counterpart in classical estimation theory.
The critical size $\norm{\Sigma}\thicksim \eps_{n,h}^2\thicksim n^{-\beta/(\beta+2)}$ justifies our specific choice
of weights $(w_j)$.
\end{example}
\begin{remark}
Our analysis extends to tests based on convex combinations of $\lambda_{r+1}(\widehat C^{t,h}),\ldots,\lambda_d(\widehat C^{t,h})$ rather than just the single eigenvalue $\lambda_{r+1}(\widehat C^{t,h})$. This approach is analogous to the likelihood ratio tests studied by \citet{onatski2014}, who consider different alternatives and their high-dimensional asymptotics. Such tests optimize power against alternatives where not only $\lambda_{r+1}(\Sigma(t))$ but also further eigenvalues $\lambda_{j}(\Sigma(t))$, $j=r+2,\ldots,d$ are significantly positive. For a fixed remaining dimension $d-r$, the asymptotic theory developed here can deliver the same detection rates for these tests.
\end{remark}

\section{Global test}
\label{SecUMTests}

We  now test $\rank(\Sigma(t))\le r$ uniformly on the entire interval $t\in [0,1]$. For this we assume integer $h^{-1}\in\N$ and split the interval into the $h^{-1}$ blocks $I_{0,h},I_{h,h}\ldots,I_{1-h,h}$. Testing ${\mathcal H}_0$ globally over $[0,1]$ is achieved by a simultaneous test on all blocks with a suitable Bonferroni correction. Our tests will have almost optimal power against local alternatives, where $\lambda_{r+1}(\Sigma(t))\ge\underline\lambda_{r+1}$ holds at least on a small interval.
By Theorem \ref{Propblockwisetest0} with level $\alpha h$ and a union bound over the $h^{-1}$ intervals $I_{kh,h}$, we obtain the following non-asymptotic result.

\begin{proposition}\label{PropUniformH0}
Suppose $h^{-1}\in\N$. For given $\alpha\in(0,1)$ the test
\begin{align*}
\phi_{\alpha}^g&:={\bf 1}\big(\exists k=0,\ldots,h^{-1}-1:\lambda_{r+1}(\widehat C^{kh,h})>\kappa_{\alpha}^g\big)\quad\text{ with}\\
\kappa_{\alpha}^g&=
  \big(4\BIAS_0^\Sigma(h)+\eps_{n,h}^2\tfrac{M^2}2\big) \Big(1+
 \log(4(d-r))\Psi\big(\tfrac{6d}{M}\big)
 +\tfrac{21\log(2(\alpha h)^{-1})}M \Big),
\end{align*}
possesses non-asymptotic level $\alpha$ uniformly over  ${\mathcal H}_0([0,1],\beta,L,\underline\lambda_r)$.
\end{proposition}

The test $\phi_{\alpha}^g$ keeps the level $\alpha$ non-asymptotically, but is clearly conservative. As a maximum of an increasing number of non-trivial blockwise tests, $\phi_{\alpha}^g$ does not seem to allow for simple asymptotically exact critical values. Note also that the convergence of extreme value statistics is already in the i.i.d. case very slow. If the bias does not dominate, we therefore recommend to use simulation-based critical values, that is the $(1-\alpha)$-quantile of $\max(\Lambda_{M,1},\ldots,\Lambda_{M,h^{-1}})$ where $\Lambda_{M,1},\ldots,\Lambda_{M,h^{-1}}$ are independent copies of $\Lambda_M$ from \eqref{eqLambdaM}. We omit the details.

By construction, $\phi_\alpha^g$ accepts ${\mathcal H}_0$ only when each blockwise test $\phi_\alpha^{(k)}$ on $I_{kh,h}$ accepts. This shows that under the global alternative
\[ {\mathcal H}_1^g(h,\underline\lambda_{r+1}):=\bigcup_{k=0}^{h^{-1}-1}{\mathcal H}_1(I_{kh,h},\underline\lambda_{r+1})
\]
the power of $\phi_\alpha^g$ is at least as large as for the blockwise tests $\phi_\alpha^{(k)}$, compare Theorem \ref{ThmH1blockwise}. Hence, inserting the asymptotics for the critical value $\kappa_\alpha^g$ and optimising the blocksize $h$, the following result is derived exactly as Corollary \ref{CorBlockTest}.

\begin{proposition}\label{PropUnif}
Consider the global test $\phi_{\alpha}^g$ for fixed $\alpha\in(0,1)$ and $M\ge r+1$. Choose the block size $h\thicksim (\underline\lambda_{r}^{-1}n)^{-1/(2\beta+2)}\vee  n^{-1/(\beta+2)}$ as $n\to\infty$. Then
 the tests $\phi_{\alpha}^g$ are uniformly consistent over local alternatives:
\[ \lim_{n\to\infty}\inf_{\Sigma\in {\mathcal H}_1^g(h,R_nv_n^g)} \PP_\Sigma(\phi_\alpha^g=1)=1
\]
for any sequence $R_n\to\infty$ with the global detection rate
\[ v_n^g=\Big(\big((\underline\lambda_{r}^{-1}L^2)^{1/(\beta+1)}n^{-\beta/(\beta+1)}\big)\wedge  \big(L^{2/(\beta+2)} n^{-\beta/(\beta+2)}\big)\Big)\log(n).\]
\end{proposition}

\begin{remark}
It is classical that the optimal detection rate under supremum norm increases by a log-factor due to multiple testing  \citep[Thm. 3.9]{ingster2012}. From classical intuition, we might want to choose $M\thicksim \log(h^{-1})\to\infty$ to balance the multiplicity correction and then equilibrate the bias with the noise level $\log(h^{-1})^2\eps_{n,h}^2$. By Remark \ref{RemH1beta} below, this is only feasible if $r=d-1$ or if we assume additionally $\Sigma\in C^\beta(L)$  under the alternative and  $\kappa_\alpha^u\gtrsim h^\beta$ (i.e., no spectral gap). In the latter case this construction yields the slightly better detection rate $(n/\log(n)^2)^{-\beta/(\beta+2)}$.
\end{remark}

\section{\bf Optimality of the detection rate}\label{SecLowerBound}

The detection rate $v_n$ for our blockwise testing problem is minimax optimal in Ingster's sense \citep{ingster2012} as the following asymptotic lower bound shows.

\begin{theorem}\label{ThmLowerBoundBlock}
Fix $t\in[0,1)$,  $\eta>0$, assume $r\ge 1$ and consider the asymptotics $n\to\infty$, $h_n\to0$.
When the local alternative ${\mathcal H}_1(I_{t,h_n},c v_n)$ is separated from the null hypothesis ${\mathcal H}_0([0,1],\beta,L,\underline\lambda_{r})$ with  detection rate $v_n$ from \eqref{Eqvn} and a small constant $c>0$, then no level-$\alpha$ test is consistent. More precisely, for any $\alpha'\in(0,1)$ there is a constant $c_{\alpha'}>0$ such that
\[\liminf_{n\to\infty}\inf_{\phi_n}\Big(\sup_{\Sigma\in {\mathcal H}_0([0,1],\beta,L,\underline\lambda_{r})}\E_{\Sigma}[\phi_n]+\sup_{\Sigma\in {\mathcal H}_1(I_{t,h_n},c_{\alpha'} v_n)}\E_{ \Sigma}[1-\phi_n]\Big)\ge\alpha',\]
 where the infimum is taken over all tests $\phi_n$ based on observing \eqref{EqY}
 and the interval length satisfies $h_n\le c_{\alpha'}^{1/(2\beta)}((\underline\lambda_r^{-1} L^2 n)^{-1/(2\beta+2)}\vee (Ln)^{-1/(\beta+2)})$.
\end{theorem}

\begin{proof}
See Appendix \ref{SecProofLB}.
\end{proof}

In essence, this theorem yields a lower bound for testing the rank of the spot covariance matrix $\Sigma(s)$ locally around a point $t$. It shows that our test has optimal asymptotic power for this problem and in particular that the smoothness $\beta$ of $\Sigma$ and the spectral gap $\underline\lambda_r$ indeed play a key role for the detection of non-zero $(r+1)$st eigenvalues.
\begin{remark}
A fine point is that Corollary \ref{CorBlockTest} yields the detection rate $v_n$ for the level-$\alpha$ test $\phi_{\alpha,0}$ under  ${\mathcal H}_0(I_{t,h},\beta,L,\underline\lambda_{r})$, but clearly $\phi_{\alpha,0}$ keeps its level under the stronger hypothesis ${\mathcal H}_0([0,1],\beta,L,\underline\lambda_{r})$. The upper bound on $h_n$ in Theorem \ref{ThmLowerBoundBlock} suggests that, for alternatives where the $(r+1)$st eigenvalue of $\Sigma(t)$ is uniformly large over a larger block, the separation rate can be improved. To exploit this non-local deviation, a test should be constructed using the average $\sum_k\lambda_{r+1}(\widehat C^{kh,h})h$  with $h$ equal to this upper bound. See \cite{RW23} for such a mean test. Yet, its analysis in our noisy observation model appears to introduce intricate technical challenges. Moreover, a mean test inherently loses power for alternatives where the $(r+1)$st eigenvalue is only locally significant---as in our data example of Section \ref{SecImpl}.
\end{remark}


Concerning the global testing problem, ${\mathcal H}_1(I_{kh_n,h_n},\underline\lambda_r)\subset {\mathcal H}_1^g(h_n,\underline\lambda_r)$ directly implies under the assumptions of Theorem \ref{ThmLowerBoundBlock}
\begin{equation}\label{EqLBuniform}
\liminf_{n\to\infty}\inf_{\phi_n}\Big(\sup_{\Sigma\in {\mathcal H}_0([0,1],\beta,L,\underline\lambda_{r})}\E_{\Sigma}[\phi_n]+\sup_{\Sigma\in {\mathcal H}_1^g(h_n,c_{\alpha'} v_n)}\E_{ \Sigma}[1-\phi_n]\Big)\ge\alpha'.
\end{equation}
We therefore lose at most the logarithmic factor of the separation rate $v_n^g$ for $\phi_\alpha^g$.

For $\underline\lambda_{r+1}=0$ (no spectral gap) under ${\mathcal H}_0$ we can consider the more canonical uniform alternative
\[ {\mathcal H}_1^{g_{can}}(\beta, L,\underline\lambda_{r+1}):=\Big\{\Sigma\in C^\beta([0,1],L)\,|\,\exists t\in[0,1]:\lambda_{r+1}(\Sigma(t))\ge\underline\lambda_{r+1}\Big\}
\]
that the $(r+1)$st eigenvalue of $\Sigma(t)$ is larger than $\underline\lambda_{r+1}$ at some point $t$ in $[0,1]$, imposing the same H\"older assumption as under ${\mathcal H}_0$. Then the inclusions
\begin{equation}\label{EqIncl}
{\mathcal H}_1^g(h_n,c v_n)\cap C_\beta([0,1],L)\subset {\mathcal H}_1^{g_{can}}(\beta,L,c v_n) \subset {\mathcal H}_1^g(h_n,(c/2) v_n)
\end{equation}
hold for any $c>0$ and $h_n=(cv_n/(2L))^{1/\beta}$ by H\"older continuity. Since we have $v_n=L^{2/(\beta+2)} n^{-\beta/(\beta+2)}$ for $\underline\lambda_r=0$, we can consider $h_n=(c/2)^{1/\beta}(Ln)^{-1/(\beta+2)}\le c_{\alpha'}^{1/(2\beta)}(Ln)^{-1/(\beta+2)}$ for a sufficiently small $c>0$. So,  for $\underline\lambda_r=0$
Proposition  \ref{PropUnif}  yields tests ${\mathcal H}_0([0,1])$ versus ${\mathcal H}_1^{g_{can}}$ at detection rate $v_n^g$. By \eqref{EqLBclasses} the lower bound in Theorem \ref{ThmLowerBoundBlock} even applies for the smaller alternative ${\mathcal H}_1\cap C^\beta([0,1],L)$. This  and \eqref{EqLBuniform} therefore show that the rate $v_n^g$ is indeed optimal for ${\mathcal H}_0([0,1],\beta,L,0)$ versus ${\mathcal H}_1^{g_{can}}(\beta,L,v_n^g)$ up to the logarithmic factor.

\begin{remark}
The testing problem is also well-defined for $r=0$, corresponding to the null hypothesis $H_0:\Sigma\equiv 0$. The lower bound proof, however, fails in this case, and standard parametric methods yield a detection lower bound of order $n^{-1}$. Surprisingly, this $n^{-1}$-rate is achieved almost trivially. By convexity, we have $\int_0^1 \lambda_1(\Sigma(t))dt\ge \lambda_1(\int_0^1\Sigma(t)dt)$, so it suffices to test $H_0:\Sigma\equiv 0$ versus  $H_1:\lambda_1(\int_0^1\Sigma(t)dt)\ge \underline\lambda_1$. That is, we only need to test the integrated covariance. This is accomplished by our test with $I_{t,h}=[0,1]$ whenever $\underline\lambda_1\gtrsim n^{-1}$, which follows from our analysis by inserting the trivial bound $\BIAS_0^\Sigma(h)=0$ and noting $\eps_{n,h}^2\thicksim n^{-1}$ for $h\thicksim 1$.
For $d=1$ and standard models it is well known that a parametric detection rate is natural \citep{ingster2012}. In contrast,  for ranks $r\ge 1$ the cone of positive semi-definite matrices $\Sigma$ with $\rank(\Sigma)\ge r+1$ is  sufficiently rich to induce our nonparametric detection rates for $r\ge 1$.
\end{remark}

\section{Implementation and examples} \label{SecImpl}

We study the finite sample performance of the asymptotic rank tests $\phi_{\alpha,1}$, $\phi_{\alpha,2}$ from Theorem \ref{Thmblockwisetest}.
In dimension $d=10$ we simulate the observations $(Y_i)$ according to \eqref{EqX} and the regression model \eqref{EqYRegr}.
The spectral statistics $S_j$ from \eqref{EqSj} are naturally approximated by
\[ \tilde S_j:=\sum_{i:(i-1/2)/n\in I_{t,h}}\Phi_j((i-1/2)/n)(Y_i-Y_{i-1}), \]
using integration by parts and Riemann sums, compare \citet{bibinger2015}. We then obtain the empirical covariance matrix $\tilde C^{t,h}=\sum_{j=1}^J\tilde w_j\tilde S_j^{\otimes 2}$ with weights $\tilde w_j\propto M^{-1}(1+j^2/M^2)^{-2}$ summing to one. Here, we choose $J=15$ and $M=10$. Variation of $J$ has little effect on the test result while with larger $M$ tests tend to be more conservative.

The sample size is $n=32,400$, referring to a one-second observation scheme over nine hours of a trading day in financial markets. Therefore the number of observations $nh$ provides the length of a block $I_{t,h}$ in seconds. The noise level is $\eta=0.001$, corresponding to the empirical findings below. We consistently test for rank $r\le 1$, i.e. ${\mathcal H}_0: \rank(\Sigma(t))\le 1$.

The simulated instantaneous covariance matrix $\Sigma(t)$ is generated by a rank-deficient Wishart process. Under ${\mathcal H}_0$ we consider $\Sigma(t)=\diag(1\;0\;\ldots\;0)+\tilde B(t)\tilde B(t)^\top$  with an independent $d$-dimensional Brownian motion $\tilde B$.
Under ${\mathcal H}_1$ we consider $r=2$ and $\tilde\Sigma(t)=\diag(1\; 0.5\; 0\;\cdots\; 0)+\tilde B(t)\tilde B(t)^\top$ with an independent $(d\times 2)$-Brownian matrix $\tilde B(t)$, containing $2d$ independent Brownian motions as entries. For different scenarios of detectability, the eigenvalues of $\tilde\Sigma(t)$ are adjusted. Given the eigen-decomposition $\tilde\Sigma(t)=V(t)\diag(\lambda_1(t),\,\lambda_2(t)) V(t)^\top$ with orthogonal $V(t)\in\mathbb R^{d\times 2}$, we define $\Sigma(t):=V(t)\diag(\lambda_1(s),\,\lambda_2^\ast)V(t)^\top$ with prescribed (non-time varying) second eigenvalue $\lambda_2^\ast>0$.

\begin{figure}[t]
\begin{center}
\includegraphics[width=0.49\textwidth]{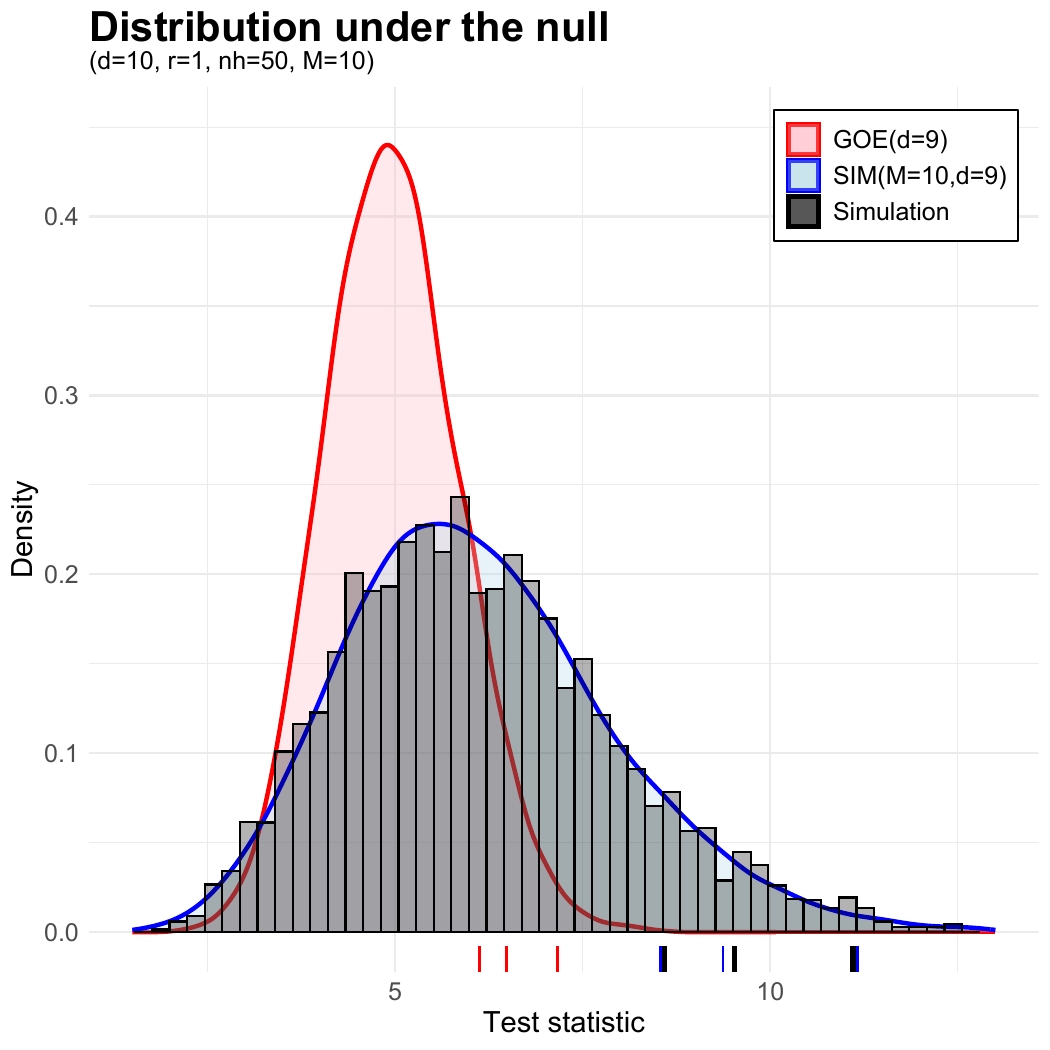}\hspace{0.1cm}
\includegraphics[width=0.49\textwidth]{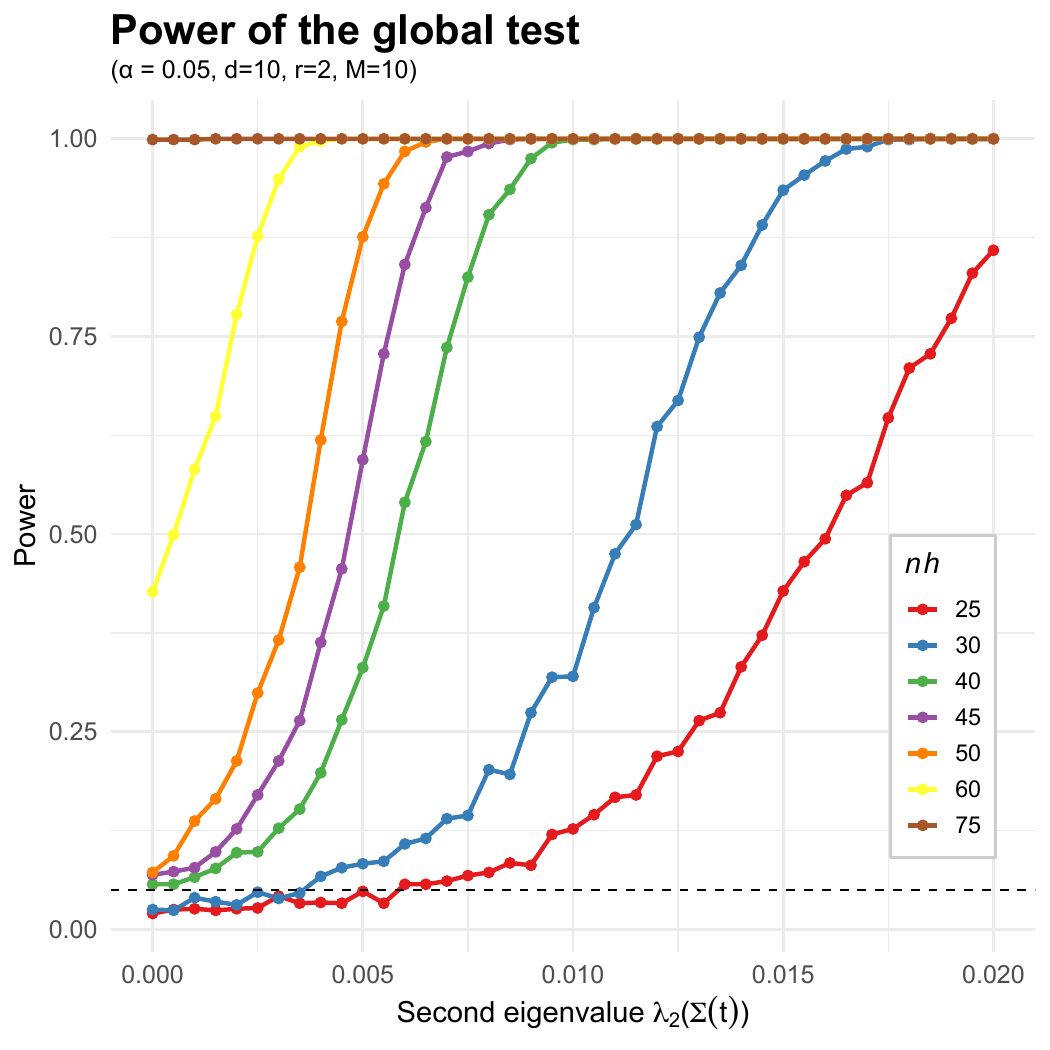}
\end{center}
\caption{ Left: Simulations under ${\mathcal H}_0$. Distribution of the largest eigenvalue of GOE in dimension 9 (red), simulation-based (SIM) distribution (blue), and the simulated second normalized eigenvalue (histogram), based on 20,000 Monte Carlo rounds. Coloured bars on the $x$-axis mark 90\%, 95\%, and 99\% quantiles. Right: rejection frequencies of the global, simulation-based test under ${\mathcal H}_1$ (1,000 repetitions).   }
\label{figh0}
\end{figure}

 The gray histogram in Figure \ref{figh0} (left) illustrates simulations under ${\mathcal H}_0$ on 50 seconds blocks with the standardised block-wise test statistic $ (\lambda_2(\widehat C^{t,h})-B_w^{(M)}\varepsilon^2_{n,h})(2\pi^{-1}M^3\varepsilon_{n,h}^4)^{-1/2}$ considered in Theorem \ref{Thmblockwisetest}.
For comparison, the blue and red densities indicate largest eigenvalue distributions of the GOE in dimension $d-r=9$ and the simulation-based test with $(\Lambda_M-B_w^{(M)})(2\pi^{-1}M^3)^{-1/2}$, respectively.

The plot suggests that while the quantiles of the simulation-based test statistics match quite well, the GOE critical values are too small, leading to elevated rejection frequencies. Rejection frequencies of the simulation-based test under ${\mathcal H}_0$ at levels $\alpha\in\{0.1; 0.05; 0.01\}$ are 10.4\%, 5.5\% and 0.9\%, respectively. The more skewed simulation-based (blue) distribution arises naturally because simulation-based critical values explicitly draw upon the sum of spectral covariance matrices across the $J=15$ ($M=10$) frequencies, while the GOE provides a limiting distribution for $J,M\to\infty$. Convergence of quantiles of eigenvalues of Wishart matrices to their Tracy-Widom analogues is  slow: for $d=10$ and even $M=1000$ tests at 10\% and 5\% level will still be oversized, see also Table 1 in \citet{ma2012}.


Our power analysis focuses on the global, simulation-based test  at $\alpha=0.05$. Figure \ref{figh0} (right) shows rejection frequencies under the rank two alternative for seven different block lengths between 25 and 75 seconds. The number of blocks on the entire interval thus lies between 1296 for $nh=25$ and 432 for $nh=75$. The plots indicate that the block length $h$ that balances $\BIAS_0^\Sigma(h)$ and $\eps_{n,h}^2$ may lay in between intervals of length 30 and 60 seconds. For $nh=25$ (red) the observation-noise part is large due to little smoothing, such that the test displays lower power. For $nh=60$ (yellow) and larger blocks the bias dominates and significantly pulls up second eigenvalues of integrated covariance matrices, even when the second eigenvalue of all spot covariance matrices is zero. This indicates that the bias assumption of Theorem \ref{ThmH1blockwise} is violated.  Block lengths between 40 and 50 seconds yield adequate size control and sensitivity.

\begin{figure}[t]
\includegraphics[width=0.49\textwidth]{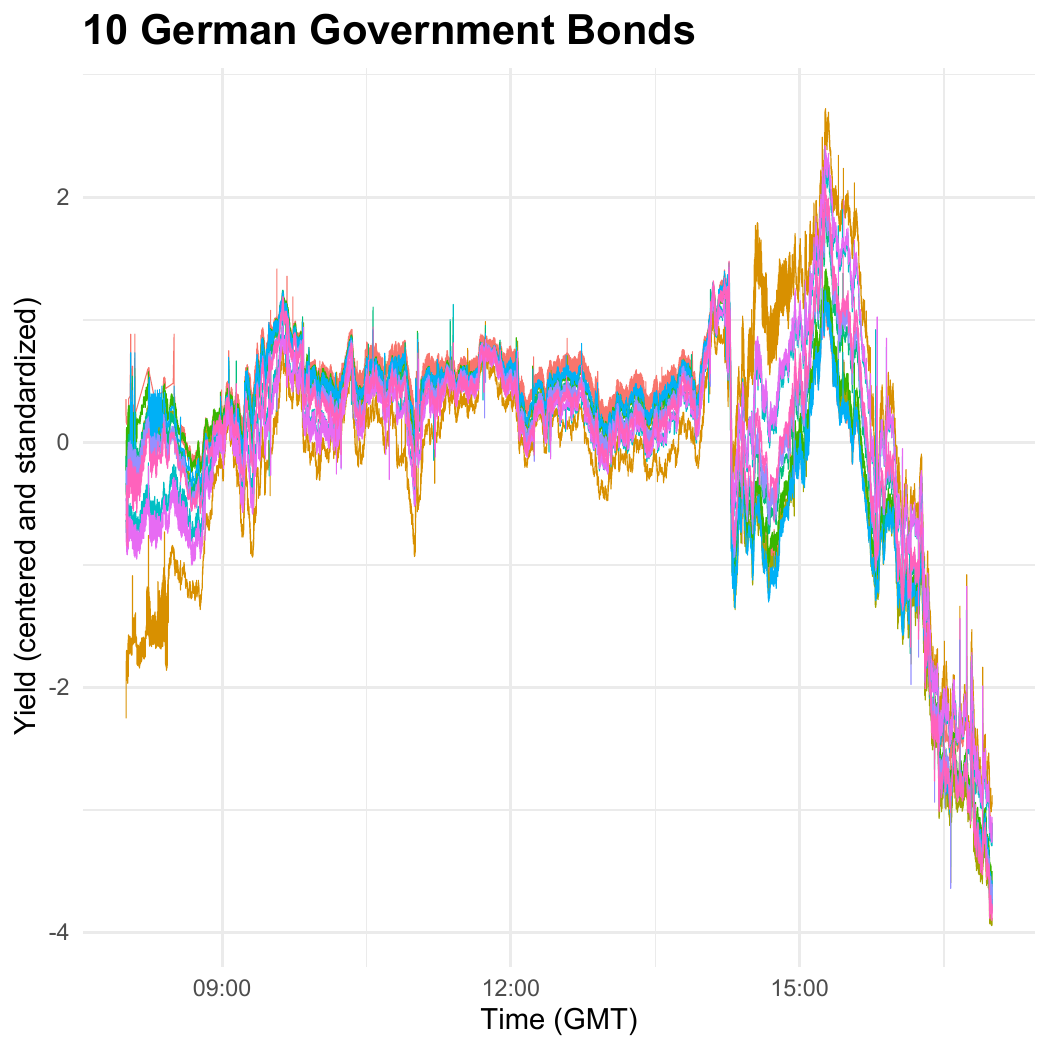}\hspace{0.1cm}
\includegraphics[width=0.49\textwidth]{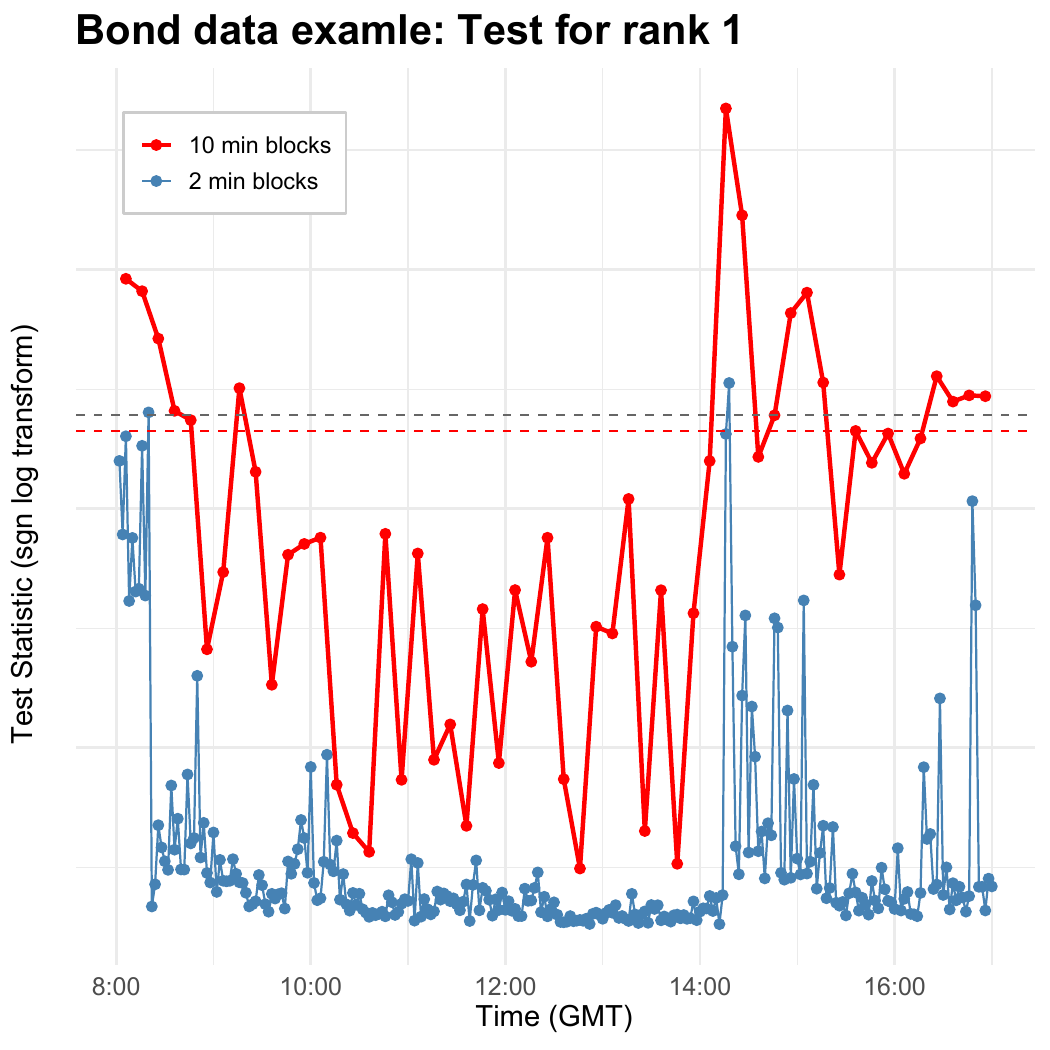}
\caption{ Data example. Left: Mid-quotes of 10 German zero coupon bonds with maturities between 3 and 30 years on 4 May 2023.  Right: corresponding test statistic  and critical values (horizontal dashed lines) of the global test.   }
\label{figh1}
\end{figure}

In a real data example we study German yield curve data, observed on 4 May 2023 (9 trading hours) and shown in Figure  \ref{figh1}  (left). Quotes are merged according to Reuters Instrument Code: primary source (PS), Frankfurt exchange (F) and composite (.). On that day the European Central Bank increased its key policy rate and indicated further rate hikes. The decision was announced at 14:15 GMT. We selected 10 zero-coupon bonds with maturities between 3 and 30 years. 
The average number of quotes for the 10 bonds on that day is 44,142. Yield changes have a first-order autocorrelation of up to -0.6, an evidence for microstructure noise. Estimated noise standard deviation is on average 0.001. For implementation details, especially in the case of asynchronous high-frequency observations, we refer to \citet{bibinger2015}.

We study covariance matrices of the 10 bonds on blocks of different length.  Figure \ref{figh1} (right) shows the test statistic for testing rank 1 using 2 and 10 minutes blocks. The dashed horizontal lines are the critical value from the  global simulation-based test. The figure shows that averaging over longer blocks increases the test statistic, leading to more frequent rejections across blocks.

\begin{table}[t]
\caption{Covariance rank for 10 bonds on  4 May 2023.}
\begin{tabular}{|lc|cccccc|}
\hline
$h$ & $K$  & rank 1& rank 2& rank 3&rank 4&rank 5&rank 6 \\ \hline
1 min & 540    &0.998  &0.002 &0&0&0&0 \\
2 min & 270    &0.993  &0.007 &0&0&0&0 \\
5 min & 108    &0.954  &0.046 &0&0&0&0 \\
10 min &54    &0.703 &0.278& 0.019&0&0&0\\
1 hour & 9 &0   &0.556 &0.222&0.111&0.111&0 \\
9 hours & 1   & 0 &0&0&0&0&1.00 \\ \hline
\multicolumn{8}{p{10cm}}{Global rank test at level $0.05$ with $K$ blocks of length $h$. Numbers give the fraction of  blocks with maximally accepted rank $r$.}
\end{tabular}

\label{tab:data}
\end{table}

Table \ref{tab:data} confirms that relation. It shows the detected rank on the specific day. A clear pattern emerges: as the block length increases, the detected rank rises from two, when using 540 one-minute blocks, to six when the block length spans the entire nine-hour trading day.
This gives significant empirical evidence that the maximal rank of the instantaneous covariance matrices can be substantially smaller than the rank of the integrated covariance matrix.

\pagebreak

\begin{appendix}

\section{Further results and proofs}

\subsection{Proof of the central limit theorem}\label{AppSecProof2}

\begin{proof}[Proof of Theorem \ref{ThmCLT}]
In order to establish asymptotic normality for (standardised) $\widehat\Sigma^{t,h}$, it suffices to establish asymptotic normality of the linear forms $\scapro{\widehat\Sigma^{t,h}}{S}_{HS}$ for all deterministic symmetric $S\in\R^{d\times d}$, in view of the Cram\'er-Wold device and the symmetry of $\widehat\Sigma^{t,h}$. Writing $S=O^\top D O$ with a diagonal matrix $D$ and an orthogonal matrix $O$, this becomes $\scapro{O\widehat\Sigma^{t,h} O^\top}{D}_{HS}$. Now observe that the estimator $\widehat\Sigma^{t,h}$ is equivariant in the sense that the law of $O\widehat\Sigma^{t,h}O^\top$ under $P_\Sigma$ is the same as the law of $\widehat\Sigma^{t,h}$ under $P_{O\Sigma O^\top}$ for any orthogonal matrix $O\in\R^{d\times d}$. This means that asymptotic normality of $\scapro{O\widehat\Sigma^{t,h} O^\top}{D}_{HS}$ under $P_\Sigma$ follows from asymptotic normality of $\scapro{\widehat\Sigma^{t,h}}{D}_{HS}$ under $P_{O\Sigma O^\top}$. It thus suffices to prove asymptotic normality of $(\widehat\Sigma^{t,h}_{\ell\ell})_{\ell=1,\ldots,d}$ under $P_{O\Sigma O^\top}$ for any orthogonal matrix $O$.

Taking into account that  $\widehat\Sigma^{t,h}$ is a quadratic form in the $S_j$ and thus an element of the second Wiener chaos generated by $B$ and $W$, a consequence of Thm. 6.2.2 in \citet{Nourdin2012} is that it suffices to establish asymptotic normality of the (individually standardised) components. We have thus reduced the problem to proving asymptotic normality of the (standardised) diagonal entries $\widehat\Sigma^{t,h}_{\ell\ell}$ under $P_\Sigma$ for any covariance function $\Sigma$ considered in the assertion of the theorem, whose conditions are orthogonally invariant.

Writing $dX(t)=\Sigma(t)^{1/2}dB(t)$ we obtain $\widehat\Sigma^{t,h}_{\ell\ell}-\E[\widehat\Sigma^{t,h}_{\ell\ell}]$ as an iterated stochastic integral:
\begin{align*}
\widehat\Sigma^{t,h}_{\ell\ell}-\Sigma^{t,h}_{\ell\ell} &= \int_{I_{t,h}}\int_{I_{t,h}}  \bscapro{\sum_{j\ge 1}w_j f_j(u)f_j(s)^\top}{d\begin{pmatrix}B(u)\\W(u)\end{pmatrix}d\begin{pmatrix}B(s)\\W(s)\end{pmatrix}^\top}_{HS}\\
\text{with } f_j(u)&:=\begin{pmatrix}\Phi_j(u)\Sigma^{1/2}(u)e_\ell\\\Phi_j'(u)\eta n^{-1/2}e_\ell\end{pmatrix}\in\R^{2d}\text{ and the $\ell$-th unit vector }e_\ell.
\end{align*}
In the parlance of Section 2.7.4 of \citet{Nourdin2012} the associated Hilbert-Schmidt operator $A_f:L^2(I_{t,h};\R^{2d})\to L^2(I_{t,h};\R^{2d})$ is  given by
\[
A_fg=\sum_{j\ge 1}w_j \scapro{f_j}{g}_{L^2(I_{t,h};\R^{2d})} f_j=\Big(\sum_{j\ge 1}w_j f_j^{\otimes 2}\Big) g.
\]
Below, we shall prove that $\norm{A_f}/\norm{A_f}_{HS}\to 0$. This analogue of the Lindeberg condition  yields that all cumulants of order three and higher for the standardisation $\Var(\widehat\Sigma^{t,h}_{\ell\ell})^{-1/2}(\widehat\Sigma^{t,h}_{\ell\ell}-\Sigma^{t,h}_{\ell\ell})$ converge to zero \cite[Prop. 2.7.13(2)]{Nourdin2012}. Consequently, asymptotic normality follows for $\widehat\Sigma^{t,h}_{\ell\ell}$  and by the above reduction for the entire matrix $\widehat\Sigma^{t,h}$.

We proceed to proving $\norm{A_f}/\norm{A_f}_{HS}\to 0$.
We calculate by considering $g=(g_1,g_2)^\top$ with $g_i\in L^2(I_{t,h};\R^d)$ and using the orthonormality of $(\Phi_j)$ and $(\Phi_j'/\norm{\Phi_j'})$, respectively,
\begin{align*}
\norm{A_f} &=\sup_{\norm{g}\le 1} \sum_{j\ge 1}w_j\scapro{f_j}{g}^2\\
&\le 2\Big(\sup_{\norm{g_1}\le 1} \sum_{j\ge 1}w_j \Big(\int_{I_{t,h}}\Phi_j(u)\scapro{\Sigma^{1/2}(u)e_\ell}{g_1(u)}\,du\Big)^2\\
&\qquad+\sup_{\norm{g_2}\le 1} \sum_{j\ge 1}w_j\Big(\eta n^{-1/2}\int_{I_{t,h}}\Phi_j'(u)\scapro{e_\ell}{g_2(u)}\,du\Big)^2\Big)\\
&\le 2 \Big(\norm{(w_j)}_{\ell^\infty}\sup_{\norm{g_1}\le 1}\int_{I_{t,h}}\scapro{\Sigma^{1/2}(u)e_\ell}{g_1(u)}^2 du+\eta^2 n^{-1} \bnorm{\big(w_j\norm{\Phi_j'}_{L^2}^2\big)}_{\ell^\infty}\Big)\\
&= 2 \big(w_1\norm{\Sigma^{1/2}e_\ell}_{L^\infty(I_{t,h})}^2 +\eps_{n,h}^2\norm{(j^2w_j)}_{\ell^\infty}\big).
\end{align*}
We now use $\norm{\Sigma^{1/2}(u)e_\ell}^2=\Sigma_{\ell,\ell}(u)$ and  inject  $w_1\le 4M^{-1}/\pi$, $\norm{(j^2w_j)}_{\ell^\infty}\le 4M/\pi$ to arrive at
\begin{equation}\label{EqNormAf}
\norm{A_f}\le \tfrac{8}{\pi}\big(M^{-1}\norm{\Sigma_{\ell,\ell}}_{L^\infty(I_{t,h})} + M\eps_{n,h}^2\big).
\end{equation}

The Hilbert-Schmidt norm of the integral operator $A_f$ is given by the $L^2(I_{t,h}\times I_{t,h};\R^{2d\times 2d})$-norm of the matrix-valued kernel \citep[Thm. 11.3.6]{BiSo2012} and $\norm{\Sigma^{1/2}(u)e_\ell}^2=\Sigma_{\ell,\ell}(u)$ holds so that
\begin{align*}
\norm{A_f}_{HS}^2 &= \int_{I_{t,h}}\int_{I_{t,h}} \bnorm{\sum_{j\ge 1} w_j f_j(u)f_j(s)^\top}_{HS(\R^{2d\times 2d})}^2duds\\
&=\sum_{j,j'\ge 1} w_jw_{j'} \int_{I_{t,h}}\int_{I_{t,h}} \trace\big( f_j(u)f_j(s)^\top f_{j'}(s)f_{j'}(u)^\top\big)\,duds\\
&=\sum_{j,j'\ge 1} w_jw_{j'} \Big(\int_{I_{t,h}} \scapro{f_j(u)}{f_{j'}(u)}\,du\Big)^2\\
&\ge \sum_{j\ge 1} w_j^2\Big(\int_{I_{t,h}}\Big(\Phi_j(u)^2\norm{\Sigma^{1/2}(u)e_\ell}^2+ \Phi_j'(u)^2\eta^2n^{-1}\Big)\,du\Big)^2\\
&= \sum_{j\ge 1} \Big(w_j\int_0^1\big(1-\cos(2\pi j s)\big)\Sigma_{\ell,\ell}(t+hs)\,ds+ \eps_{n,h}^2w_jj^2\Big)^2,
\end{align*}
compare \eqref{EqSigmath}. Now observe that by the orthonormality of the cosine functions
\[ \sum_{j\ge 1} \Big(\int_0^1\cos(2\pi j s)\Sigma_{\ell,\ell}(t+hs)\,ds\Big)^2\le \norm{\Sigma_{\ell,\ell}(t+h\cdot)}_{L^2([0,1])}^2\]
holds. By the inverse triangle inequality in $\ell^2$ we can thus lower bound
\begin{align*}
&\Big(\sum_{j\ge 1} \Big(w_j\int_0^1\big(1-\cos(2\pi j s)\big)\Sigma_{\ell,\ell}(t+hs)\,ds\Big)^2\Big)^{1/2}\\
&\ge \norm{(w_j)}_{\ell^2}
\norm{\Sigma_{\ell,\ell}(t+h\cdot)}_{L^1([0,1])}-\norm{(w_j)}_{\ell^\infty}
\norm{\Sigma_{\ell,\ell}(t+h\cdot)}_{L^2([0,1])}.
\end{align*}
Using $\norm{(w_j)}_{\ell^2}\thicksim M^{-1/2}$, $\norm{(w_j)}_{\ell^\infty}\thicksim M^{-1}$, $\norm{(j^2w_j)}_{\ell^2}\thicksim M^{3/2}$, we conclude
\begin{align}
\norm{A_f}_{HS}&\gtrsim \max\Big(M^{3/2}\eps_{n,h}^2,\label{EqHSNormAf}\\
&\qquad M^{-1/2}\norm{\Sigma_{\ell,\ell}(t+h\cdot)}_{L^1([0,1])}-CM^{-1}\norm{\Sigma_{\ell,\ell}(t+h\cdot)}_{L^2([0,1])}
\Big)\notag
\end{align}
for some constant $C>0$.

From the condition \eqref{EqAssCLT} with $v=e_\ell$ we have $\norm{\Sigma_{\ell,\ell}}_{L^\infty(I_{t,h})}\lesssim M^2\eps_{n,h}^2$ or, using $\norm{\Sigma_{\ell,\ell}(t+h\cdot)}_{L^p([0,1])}=h^{-1/p}\norm{\Sigma_{\ell,\ell}}_{L^p(I_{t,h})}$ for $p\in[1,\infty]$ and $L^p$-interpolation,
\[ \norm{\Sigma_{\ell,\ell}(t+h\cdot)}_{L^\infty([0,1])}\lesssim \norm{\Sigma_{\ell,\ell}(t+h\cdot)}_{L^2([0,1])}\lesssim \norm{\Sigma_{\ell,\ell}(t+h\cdot)}_{L^1([0,1])}.
\]
We can thus deduce from \eqref{EqNormAf} and \eqref{EqHSNormAf}  \[\frac{\norm{A_f}}{\norm{A_f}_{HS}}\lesssim \frac{M^{-1}\norm{\Sigma_{\ell,\ell}}_{L^\infty(I_{t,h})}
+M\eps_{n,h}^2} {M^{-1/2}\norm{\Sigma_{\ell,\ell}}_{L^\infty(I_{t,h})}+M^{3/2}\eps_{n,h}^2}\thicksim M^{-1/2}\to 0,\]
as required.

Finally, let us establish the bound on the covariance matrix for $\vek(\widehat\Sigma^{t,h})$. To do so, it suffices, again by the equivariance with respect to orthogonal transformations, to show for any diagonal matrix $H=\sum_{\ell=1}^d h_\ell e_\ell e_\ell^\top$ that
\begin{equation}\label{EqCovBound}
\Var(\scapro{\widehat\Sigma^{t,h}}{H}_{HS})\lesssim M^{-1}\bnorm{\big((\norm{\Sigma_{\ell,\ell'}}_{L^\infty(I_{t,h})})_{\ell,\ell'}+\eps_{n,h}^2M^2I_d\big)H}_{HS}^2.
\end{equation}
Using $\Var(\scapro{\widehat\Sigma^{t,h}}{H}_{HS})\le d \sum_{\ell=1}^d h_\ell^2 \Var(\widehat\Sigma_{\ell\ell}^{kh})$ and evaluating the Hilbert-Schmidt norm on the right only on the diagonal, \eqref{EqCovBound} follows from
\begin{equation}\label{EqCovBound2}
\Var(\widehat\Sigma_{\ell,\ell}^{t,h})\lesssim M^{-1}\big(\norm{\Sigma_{\ell,\ell}}_{L^\infty(I_{t,h})}+\eps_{n,h}^2M^2\big)^2
\end{equation}
for $\ell=1,\ldots,d$.
Now, from the previous calculations we directly get
\begin{align*}
\Var(\widehat\Sigma^{t,h}_{\ell\ell}) &= 2\norm{A_f}_{HS}^2\\
&\le 2\trace(A_f)\norm{A_f}\\
&\lesssim \Big(\Sigma_{\ell\ell}^{t,h}+\eps_{n,h}^2M^2\Big)M^{-1}\Big(\norm{\Sigma_{\ell\ell}}_{L^\infty(I_{t,h})}+\eps_{n,h}^2M^2\Big)\\
&\lesssim M^{-1}\Big(\norm{\Sigma_{\ell\ell}}_{L^\infty(I_{t,h})}+\eps_{n,h}^2M^2\Big)^2,
\end{align*}
hence \eqref{EqCovBound2}.
\end{proof}

\subsection{Eigenvalue concentration}\label{secMatrixConc}

In the following we shall work on general real separable Hilbert spaces $H$. In a first step the expected spectral norm of Gaussian operators on $H$ is bounded tightly. Then the spectral norm for certain matrix-valued martingales in the second Wiener chaos is bounded via a quadratic variation argument involving a Gaussian operator. Finally, Gaussian concentration for the square-root spectral norm of matrices in the second Wiener chaos yields a  deviation bound.

The original Lieb Theorem \citep{lieb1973} considers selfadjoint operators whose exponentials are trace-class. In the next result we get around this restrictive assumption by working with the function $e^x-1-x$ instead of the exponential.

\begin{lemma}\label{LemLieb}
Consider a deterministic Hilbert-Schmidt operator $F:H\to H$ and a random Hilbert-Schmidt operator $X:H\to H$ with $\E[X]=0$ and nuclear norm $\norm{\log(\E[e^X])}_{nuc}<\infty$. Then
\[ \E[\trace(\exp(F+X)-\Id-(F+X))]\le \trace(\exp(F+\log(\E[e^X]))-\Id-F)\]
and the right-hand side is finite. Here $\Id$ denotes the identity.
\end{lemma}

\begin{proof}
On a finite-dimensional Hilbert space the assertion boils down to an application of Jensen's inequality in Lieb's Theorem \cite[Cor 3.4.2]{tropp2012}. Therefore assume $\dim(H)=\infty$ and consider an orthonormal basis $(e_k)_{k\ge 1}$ of $H$. Introduce $\trace_n(T):=\sum_{k=1}^n\scapro{Te_k}{e_k}$, the trace of the operator $T$ restricted to $\spann(e_1,\ldots,e_n)$. We have $e^T-\Id-T\ge 0$ for any self-adjoint operator $T$ because of $e^x-1-x\ge 0$ for all $x$. This implies $\trace_n(e^T-\Id-T)\uparrow \trace(e^T-\Id-T)$ as $n\to\infty$. Monotone convergence and the finite-dimensional result on $\spann(e_1,\ldots,e_n)$ therefore yield
\begin{align*}
&\E\Big[\trace\big(\exp(F+X)-\Id-(F+X)\big)\Big]\\
&= \lim_{n\to\infty} \E\Big[\trace_n\big(\exp(F+X)-\Id-(F+X)\big)\Big]\\
&\le \liminf_{n\to\infty}\trace_n\big(\exp(F+\log(\E[e^X]))-\Id-F\big)\\
&= \trace\big(\exp(F+\log(\E[e^X]))-\Id-(F+\log(\E[e^X]))\big)+\lim_{n\to\infty} \trace_n(\log(\E[e^X]))\\
&= \trace\big(\exp(F+\log(\E[e^X]))-\Id-F\big),
\end{align*}
where we used $\trace_n(T)\to \trace(T)$ for nuclear $T=\log(\E[e^X])$ in the last step. Finally, observe that the two last lines are finite because $\log(\E[e^X])$ is trace-class and $\trace(e^T-\Id-T)\le \trace(T^2)e^{\norm{T}}<\infty$ for the Hilbert-Schmidt operator $T=F+\log(\E[e^X])$, a consequence of $e^x-1-x\le x^2e^{\abs{x}}$.
\end{proof}

The following result follows similarly to the intrinsic matrix Bernstein bound in \cite[Thm. 7.3.1]{tropp2012}.

\begin{proposition}\label{PropHSIneqSA}
Let $S_k: H\to H$, $k=1,\ldots,K$, be deterministic selfadjoint Hilbert-Schmidt operators. Consider the random operator $Z=\sum_{k=1}^K\gamma_kS_k$ with independent $\gamma_k\sim N(0,1)$. For a selfadjoint trace-class operator $V$ on $H$ with $V\ge \E[Z^2]=\sum_{k=1}^K S_k^2$ let
$\rho=\trace(V)/\norm{V}$ be its intrinsic dimension.
Then for all $\tau\ge 0$
\[ \PP\Big(\lambda_{max}(Z)\ge \tau\Big)\le \frac{7\rho}{4}\exp\Big(-\frac{\tau^2}{2\norm{V}}\Big)\text{ and } \E\Big[\lambda_{max}(Z)^2\Big]\le 2\log\big(7e \rho/4\big) \norm{V}.\]
\end{proposition}

\begin{proof}
 As for symmetric matrices \cite[Lemma 4.6.2]{tropp2012} we obtain for bounded selfadjoint linear operators $S_k$
\[ \E[\exp(\theta\gamma_kS_k)]=\exp(\theta^2 S_k^2/2),\quad \theta\in\R.\]
Applying Lemma \ref{LemLieb} with $F=0$ and $X=\sum_{k=1}^K \gamma_kS_k$, using that the $\gamma_k$ are independent, we find
\begin{align*}
&\E\Big[\trace\Big(\exp\Big(\theta\sum_{k=1}^K \gamma_kS_k\Big)-\theta\sum_{k=1}^K \gamma_kS_k-\Id\Big)\Big]\\
& \le \trace\Big(\exp\Big(\sum_{k=1}^K \log\Big(\E\Big[e^{\theta \gamma_kS_k}\Big]\Big)\Big)-\Id\Big) = \trace\Big(\exp\Big(\frac{\theta^2}{2}\sum_{k=1}^K S_k^2\Big)-\Id\Big).
\end{align*}
By convexity of $x\mapsto e^x-1$ and $V\ge \sum_{k=1}^K S_k^2\ge 0$ we obtain for $\theta\ge 0$
\begin{align*}
\trace\Big(\exp\Big(\frac{\theta^2}{2}\sum_{k=1}^K S_k^2\Big)-\Id\Big)&\le \trace\Big(\exp\Big(\frac{\theta^2}{2}V\Big)-\Id\Big)\\
 &\le \trace(V)\Big(\frac{\exp(\frac{\theta^2}{2}\norm{V})-1}{\norm{V}}\Big)\le \rho\exp\Big(\frac{\theta^2}{2}\norm{V}\Big),
\end{align*}
compare \cite[Eq. (7.7.3)]{tropp2012}. For the deviation probability the approach for Thm. 7.3.1 in \citet{tropp2012} yields for $\tau>0$
\begin{align*}
\PP\Big(\lambda_{max}(Z)\ge \tau\Big) &\le \inf_{\theta> 0} (e^{\theta \tau}-\theta \tau-1)^{-1}\E[\trace(\exp(\theta Z)-\theta Z-\Id)]\\
 &\le \inf_{\theta> 0} (e^{\theta \tau}-\theta \tau-1)^{-1}\rho\exp\Big(\frac{\theta^2}{2}\norm{V}\Big)\\
 &\le \rho \inf_{\theta> 0} \Big(1+\frac{3}{\theta^2\tau^2}\Big)\exp\Big(-\theta \tau+\frac{\theta^2}{2}\norm{V}\Big)\\
  &\le \rho \Big(1+\frac{3v^2}{\tau^4}\Big)\exp\Big(-\frac{\tau^2}{2\norm{V}}\Big),
 \end{align*}
inserting $\theta=\tau/\norm{V}$. For $\tau\ge (2\norm{V})^{1/2}$ the prefactor is smaller than $7\rho/4$. This implies the deviation bound because for $\tau\le (2\norm{V})^{1/2}$ the bound $(7\rho/4)e^{-\tau^2/(2\norm{V})}$ is always larger than the trivial bound $1$.

For the second moment we bound for any $A>0$:
\begin{align*}
\E\Big[\lambda_{max}(Z)^2\Big] &\le A^2+\int_A^\infty 2\tau \PP(\lambda_{max}(Z)\ge \tau)\,d\tau\\
&\le A^2+\frac{7\rho}{2} \int_A^\infty \tau \exp\Big(-\frac{\tau^2}{2v}\Big)\,d\tau
=A^2+\frac{7\rho v}{2} \exp\Big(-\frac{A^2}{2v}\Big).
\end{align*}
For $A^2=2\log(7\rho/4)v$ the last expression equals $2\log(7e\rho/4)v$.
\end{proof}

\begin{proposition}\label{PropHSIneqGeneral}
Let $Z$ be a centred Gaussian random element in $HS(H_1,H_2)$, the Hilbert space of Hilbert-Schmidt operators between real separable Hilbert spaces $H_1,H_2$. Setting
\[v:=\max(\norm{\E[Z^\ast Z]},\norm{\E[ZZ^\ast]}), \quad \rho:=\frac{2\E[\trace(Z^\ast Z)]}{v},\]
we obtain for all $\tau\ge 0$
\[ \PP\Big(\norm{Z}\ge \tau\Big)\le \frac{7\rho}{4}\exp\Big(-\frac{\tau^2}{2v}\Big)\text{ and } \E\big[\norm{Z}^2\big]\le 2\log\big(7e\rho/4\big) v.\]
\end{proposition}

\begin{proof}
Define the dilation $\tilde Z\in HS(H_1\times H_2)$ via $\tilde Z(v_1,v_2):=(Zv_2,Z^\ast v_1)$, which is selfadjoint with $\lambda_{max}(\tilde Z)=\norm{Z}$ and $\E[\tilde Z^2](v_1,v_2)=(\E[ZZ^\ast] v_1,\E[Z^\ast Z]v_2)$, $\norm{E[\tilde Z^2]}=v$  \cite[Sec. 2.1.16]{tropp2012}.

We use the Karhunen-Lo\`eve decomposition $\tilde Z=\sum_{k=1}^\infty \lambda_k^{1/2}\gamma_k E_k$ with deterministic and orthonormal $E_k\in HS(H_1\times H_2)$, independent $\gamma_k\sim N(0,1)$ and the eigenvalues $\lambda_k$ of the trace-class covariance operator of $\tilde Z$. The approximations $\tilde Z_n:=\sum_{k=1}^n \lambda_k^{1/2}\gamma_k E_k$ have exactly the form of $Z$ in Proposition \ref{PropHSIneqSA} with $S_k=\lambda_k^{1/2}E_k$ and $\E[\tilde Z_n^2]=\sum_{k=1}^n \lambda_k E_k^2\le \E[\tilde Z^2]=:V$. The intrinsic dimension of $V$ is given by $\rho$ because of $\trace(\tilde Z)=2\trace(Z^\ast Z)$. Proposition \ref{PropHSIneqSA} therefore yields for any $n\in\N$
\begin{align*}
\PP\Big(\lambda_{max}(\tilde Z_n)\ge \tau\Big)\le \frac{7\rho}{4}\exp\Big(-\frac{\tau^2}{2v}\Big),\quad
\E\Big[\lambda_{max}(\tilde Z_n)^2\Big]\le 2\log\big(7e\rho/4\big)v.
\end{align*}
The convergence $\tilde Z_n\xrightarrow{L^2}\tilde Z$ and the continuity of $S\mapsto \lambda_{max}(S)$ on $HS(H_1\times H_2)$  entail that these bounds remain true for the limit $\lambda_{max}(\tilde Z)$. As a last step we insert the identity $\lambda_{max}(\tilde Z)=\norm{Z}$.
\end{proof}

Next, the expected spectral norm of matrix-valued continuous martingales is bounded.

\begin{proposition}\label{PropSemimartBernstein}
For $G_j\in L^2([0,1];\R^{d\times d})$, $j\ge 1$, consider
\[ Q=\sum_{j\ge 1}\int_0^1 X_j(s)\,dX_j(s)^\top\in\R^{d\times d}\text{ with } X_j(s)=\int_0^sG_j(u)\,dB(u)\in\R^d,\]
which is a well-defined random matrix if $v_1$ below is finite. Then
\begin{align*}
 \E\big[\lambda_{max}(Q+Q^\top)\big]\le \log(\sqrt{7e/2}d)\max(v_1,v_2)
 \end{align*}
holds with
\begin{align*}
v_1^2 &=\lambda_{max}\big(\E\big[(Q+Q^\top)^2\big]\big)\\
&=\bnorm{\sum_{j,j'\ge 1}\Big(\scapro{G_j}{G_{j'}}_{L^2([0,1];\R^{d\times d})}\int_0^1 G_jG_{j'}^\top+\Big(\int_0^1 G_jG_{j'}^\top\Big)^2\Big)},\\
v_2^2 &=\sup_{\norm{F}_{L^2([0,1];\R^{d\times d})}\le 1}\bnorm{\sum_{j\ge 1}\int_\cdot^1\big(\scapro{G_j(s)}{F(s)}_{HS}I_d+G_j(s)F(s)^\top\big)\, ds\,G_j}_{L^2([0,1];\R^{d\times d})}^2.
\end{align*}
 \end{proposition}

\begin{remark}
The idea of the proof is to interpret $Q+Q^\top$ as a matrix-valued martingale, to bound its maximal eigenvalue by a corresponding quadratic variation expression and to interpret this expression as the squared norm of a Gaussian linear operator $T:\R^d\to L^2([0,1];\R^{d\times d})$.  Then Proposition \ref{PropHSIneqGeneral} is applied to bound $\E[\norm{T}^2]$.
\end{remark}

\begin{proof}
We consider the symmetric matrix-valued continuous martingale
\begin{align*}
M(t) &:=\sum_{j\ge 1}\int_0^t \Big(X_j(s)dX_j(s)^\top+\big(X_j(s)dX_j(s)^\top\big)^\top \Big),\quad t\in[0,1].
\end{align*}
The $i$th column  of $M(t)$ is
$M(t)_{\cdot, i}=\sum_{j\ge 1}\int_0^t (X_j(s)\, dX_{j,i}(s)+X_{j,i}(s) dX_j(s))$
and thus the sum over its quadratic covariation matrices is
\begin{align*}
&A(t):=\sum_{i=1}^d\langle M_{\cdot,i}\rangle_t \\
&=\sum_{j,j'\ge 1}\int_0^t \Big(  X_j(s)\trace\big(G_j(s)^\top G_{j'}(s)\big)  X_{j'}(s)^\top
+X_j(s) G_j(s)^\top   X_{j'}(s)^\top G_{j'}(s) \\&\quad +G_j(s)   G_{j'}(s)^\top X_j(s) X_{j'}(s)^\top
+G_j(s)X_j(s)^\top  X_{j'}(s) G_{j'}(s)^\top\Big)\,ds.
\end{align*}
Collecting terms yields that $w^\top A(1) w$ equals for $w\in\R^d$
\begin{align*}
\int_0^1\bnorm{\sum_{j\ge 1}\Big(w^\top X_j(s)G_j(s)+X_j(s)w^\top G_j(s)\Big)}_{HS}^2ds
= \norm{Tw}_{L^2([0,1];\R^{d\times d})}^2
\end{align*}
with the Gaussian linear operator $T:\R^d\to {L^2([0,1];\R^{d\times d})}$ given by
\[ (Tw)(s):=\sum_{j\ge 1}\Big(\scapro{w}{X_j(s)} G_j(s)+X_j(s) w^\top G_j(s)\Big).\]

The It\^o formula for the trace exponential \cite[Lemma 3]{bacry2018} of a continuous It\^o-semimartingale $Z(t)\in\R^{d\times d}$ yields
\[ d(\trace(\exp(Z(t))))=\trace(\exp(Z(t))dZ(t))+\frac12\sum_{i=1}^d\trace(\exp(Z(t))d\langle Z_{\cdot,i}\rangle_t).\]
 Therefore, if we consider $Z(t)=\alpha M(t)-\frac{\alpha^2}2A(t)$ for $\alpha\in\R$,
we obtain that $\trace(\exp( Z(t)))$ forms a super-martingale with
\[\E\Big[\trace\big(\exp( Z(1))\big)\Big]\le \E\Big[\trace\big(\exp(Z(0))\big)\Big]=d.\]
We deduce with Jensen's inequality
$\E[\lambda_{max}(Z(1))]\le \log(d)$ and thus $\E[\lambda_{max}(\alpha M(1)-\tfrac{\alpha^2}{2}A(1))]\le \log(d)$.
Using the convexity of $\lambda_{max}$, this yields
\begin{align}\label{EqNormM1}
\E[\lambda_{max}(M(1))]& \le \inf_{\alpha> 0}\Big(\tfrac{\alpha}2\E[\lambda_{max}(A(1))]+\tfrac1\alpha\log d\Big)=\Big(\tfrac{\log(d)}{2}\E[\lambda_{max}(A(1))]\Big)^{1/2}.
\end{align}
With the representation of $A(1)$ via $T$ we obtain
\[ \E[\lambda_{max}(M(1))]\le \Big(\tfrac{\log(d)}{2}\E[ \norm{T}_{\R^d\to L^2([0,1];\R^{d\times d})}^2]\Big)^{1/2}.\]

From Proposition \ref{PropHSIneqGeneral} we deduce
\begin{equation}\label{EqnormT} \E\big[\norm{T}_{\R^d\to L^2([0,1];\R^{d\times d})}^2\big]\le 2\log\big(7e\rho/4\big) v
\end{equation}
with $\rho=2d$ due to  $\trace(\E[T^\ast T])\le d v$, setting
\[ v=\max\Big(\norm{\E[A(1)]},\sup_{\norm{F}_{L^2}=1}\E\big[\scapro{T^\ast F}{T^\ast F}_{\R^d}\big]\Big).
\]
The It\^o isometry yields
\begin{align*}
\E\big[A(1)\big]
&=  \sum_{j,j'\ge 1}\Big( \trace\Big(\int_0^1G_j(s)G_{j'}(s)^\top ds\Big) \int_0^1G_j(s)G_{j'}(s)^\top ds \\
&\qquad\qquad
+ \Big(\int_0^1G_{j}(s)G_{j'}(s)^\top ds\Big)^2 \Big).
\end{align*}
This gives $\norm{\E[A(1)]}=v_1^2$, where the representation as an expectation follows directly from the quadratic variation property of $A(1)$. Moreover,
$M(1)=Q+Q^\top$ is well-defined if $v_1<\infty$.
We obtain the adjoint $T^\ast$ applied to $F\in L^2([0,1];\R^{d\times d})$ as
\begin{align*}
T^\ast F&=\sum_{j\ge 1}\int_0^1\Big(\big(\scapro{G_j(s)}{F(s)}_{HS}I_d+G_j(s)F(s)^\top\big) X_j(s)\Big)\, ds\\
&=\sum_{j\ge 1}\int_0^1\int_u^1\Big(\big(\scapro{G_j(s)}{F(s)}_{HS}I_d+G_j(s)F(s)^\top\big)  \, ds\,dX_j(u).
\end{align*}
This implies by It\^o isometry
\begin{align*}
&\E\big[\norm{T^\ast F}_{\R^d}^2\big] = \int_0^1\bnorm{\sum_{j\ge 1}\int_u^1\big(\scapro{G_j(s)}{F(s)}_{HS}I_d+G_j(s)F(s)^\top\big)\, ds\,G_j(u)}_{HS}^2du.
\end{align*}
We conclude $v=\max(v_1^2,v_2^2)$ and obtain the result by inserting \eqref{EqnormT} into \eqref{EqNormM1}.
\end{proof}

\begin{theorem}\label{ThmMatrixDev}
In the notation of Proposition  \ref{PropSemimartBernstein} we have with probability at least $1-\alpha$
\begin{align*}
\lambda_{max}\Big(\sum_{j\ge 1}X_j(1)^{\otimes 2}\Big)
& \le \inf_{\delta>0}\Big((1+\delta)\Big(\lambda_{max}\Big(\sum_{j\ge 1}\E\big[X_j(1)^{\otimes 2}\big]  \Big)\\
&\qquad + \log(\sqrt{7e/2}d)\max(v_1,v_2)\Big) +(1+\delta^{-1})2\log(\alpha^{-1}) \sigma^2\Big)
\end{align*}
with
\begin{align*}
\sigma^2&= \sup_{\norm{f}_{L^2([0,1];\R^d)}\le 1}\bnorm{\sum_{j\ge 1}\Big(\int_0^1 G_jf\Big)^{\otimes 2}}.
\end{align*}
\end{theorem}

\begin{proof}
Set $Z:=\sum_{j\ge 1}X_j(1)^{\otimes 2}$ and note by It\^o's formula
\[Z-\E[Z]=\sum_{j\ge 1}\int_0^1\Big(X_j(s)dX_j(s)^\top +\big(X_j(s)dX_j(s)^\top \big)^\top\Big).
\]
Therefore, by Proposition \ref{PropSemimartBernstein} we find
\begin{align*}
\E[\lambda_{max}(Z)]&\le \lambda_{max}(\E[Z])+\log(\sqrt{7e/2}d)\max(v_1,v_2) .
 \end{align*}
We observe $Z\ge 0$ so that $\lambda_{max}(Z)\ge 0$ and we use Gaussian concentration under Lipschitz maps \cite[Eq. (1.6)]{LedouxTal1991}:
\[ \PP\Big(\lambda_{max}(Z)^{1/2}\ge \E[\lambda_{max}(Z)^{1/2}]+\sqrt{2\log(\alpha^{-1})}\sigma\Big)\le \alpha\]
with the  Lipschitz constant
\begin{align*}
\sigma=\sup_{\norm{f}_{L^2([0,1];\R^d)}\le 1}\lambda_{max}\Big(\sum_{j\ge 1}\Big(\int_0^1 G_jf\Big)^{\otimes 2}\Big)^{1/2}.
\end{align*}
We note $\E[\lambda_{max}(Z)^{1/2}]\le \E[\lambda_{max}(Z)]^{1/2}$ and obtain with probability at least $1-\alpha$
\begin{align*}
\lambda_{max}(Z)^{1/2}\le \Big(\lambda_{max}(\E[Z])+\log(\sqrt{7e/2}d)\max(v_1,v_2)\Big)^{1/2} +\sqrt{2\log(\alpha^{-1})} \sigma.
\end{align*}
It remains to square both sides, using $(A+B)^2\le(1+\delta)A^2+(1+\delta^{-1})B^2$.
\end{proof}

\subsection{Empirical eigenvalue concentration under ${\mathcal H}_0$}\label{secEVH0}

The preceding results are now applied to bound the type $I$ error. Recall  ${\mathcal W}(x)=x\vee\sqrt{x}$ for $x\ge 0$.

\begin{proposition}\label{PropVtilde}
Consider
\[ \tilde V_{>r}^{t,h}=\sum_{j\ge 1} w_j\Big(\int_{I_{t,h}} \Phi_j'(s)\tfrac{\eta}{n^{1/2}}\,dW_{>r}(s)\Big)^{\otimes 2}.\]
Then for $\alpha\in(0,1)$ we have with probability at least $1-\alpha$
\[\lambda_{max}(\tilde V_{>r}^{t,h})\le  \eps_{n,h}^2\tfrac{M^2}4\Big(1+
 \log(4(d-r)){\mathcal W}\big(6(d-r+1)M^{-1}\big)
 +21M^{-1}\log(\alpha^{-1}) \Big).
\]
\end{proposition}

\begin{proof}
It suffices to remark that the $(\Phi_j')$ are $L^2$-orthogonal and $\frac{\eta^2}n\norm{\Phi_j'}_{L^2}^2=\eps_{n,h}^2j^2$ holds so that $\tilde V_{>r}^{t,h}$ is distributed as $\sum_{j\ge 1}\Gamma_j^{\otimes 2}$ in \eqref{EqDevIneqEx} with $s_j=\eps_{n,h}^2w_jj^2$. Using $\norm{(j^2w_j)}_{\ell^1}\le M^2/8$, $\norm{(j^2w_j)}_{\ell^2}\le M^{3/2}/\sqrt{2\pi}$, $\norm{(j^2w_j)}_{\ell^\infty}\le 4M/\pi$ and $\delta=1$ in \eqref{EqDevIneqEx} yields the result when simplifying the numerical constants.
\end{proof}

\begin{proposition}\label{PropSigmatilde}
Under ${\mathcal H}_0$ consider
\[ \tilde\Sigma_{>r}^{t,h}=\sum_{j\ge 1} w_j\Big(\int_{I_{t,h}} \Phi_j(s)\Sigma_{>r}^{1/2}(s)\,dB_{>r}(s)\Big)^{\otimes 2}.\]
Then for $\alpha\in(0,1)$, $M\ge r$ we have with probability at least $1-\alpha$
\[\lambda_{max}(\tilde\Sigma_{>r}^{t,h})
  \le 2\BIAS_0^\Sigma(h)\Big(1+ 2r^{1/2}\log(4(d-r))M^{-1/2}+3\log(\alpha^{-1})M^{-1}\Big).
\]
\end{proposition}

\begin{proof}
To ease notation, assume $t=0$, i.e. $I_{t,h}=[0,h]$.
We apply Theorem \ref{ThmMatrixDev} with
\[G_j(s)=w_j^{1/2} \Sigma(s)_{>r}^{1/2}\Phi_j(s) \in\R^{(d-r)\times (d-r)},\;\Phi_j(s)=\sqrt {2/h}\sin(\pi jt/h){\bf 1}_{[0,h]}(s).
\]

For $v_1^2$ in Theorem \ref{ThmMatrixDev} we note
$\iint (\sum_{j\ge 1}w_j\Phi_j(s)\Phi_{j}(u))^2dsdu=\sum_{j\ge 1}w_j^2$
because $(\Phi_j)$ forms an $L^2$-orthonormal system and  bound
\begin{align*}
v_1^2 &=  \bnorm{\sum_{j,j'\ge 1} w_jw_{j'}\Big(\int\big(\trace_{>r}(\Sigma)\Phi_j\Phi_{j'}\big)\int \Sigma_{>r}\Phi_j\Phi_{j'}  + \Big(\int \Sigma_{>r}\Phi_j\Phi_{j'}\Big)^2\Big)}\\
&= \bnorm{\iint \Big(\sum_{j\ge 1}w_j\Phi_j(s)\Phi_{j}(u)\Big)^2\Big(\trace_{>r}(\Sigma(s))\Sigma_{>r}(u)
+\Sigma_{>r}(s)\Sigma_{>r}(u)\Big)\, dsdu}\\
&\le (r+1)\norm{\Sigma_{>r}}_{L^\infty([0,h])}^2\norm{(w_j)}_{\ell^2}^2,
\end{align*}
where we used that $\Sigma(s)$ and thus also $\Sigma_{>r}(s)$ has at most rank $r$ under ${\mathcal H}_0$.

Furthermore, we extend \eqref{Defwj} and set $ w_{j}=c_w M^{-1}(1+j^2/M^2)^{-2}$ for all $j\in\Z$. Then the sequence $(w_j)_{j\in\Z}$ is positive definite because $f(x)=(1+x^2)^{-2}$ is a positive-definite function as the characteristic function of a symmetric bilateral $\Gamma$-distribution, compare \cite{kuchler2008}. So, by Herglotz's Theorem \citep{katznelson2004} $\frac{w_0}{2}+\sum_{j\ge 1} w_j\cos(x)\ge 0$ holds for all $x\in\R$. With the orthonormal system $\Psi_j(s)=\sqrt {2/h}\cos(\pi js/h){\bf 1}_{[0,h]}(s)$ for $j\ge 1$, $\Psi_0(s)=\sqrt{1/h}{\bf 1}_{[0,h]}(s)$  we obtain for $s,u\in [0,h]$
\begin{align*}
\babs{\sum_{j\ge 1}w_j \Phi_j(s)\Phi_j(u)} &=\babs{\frac{2}{h}\sum_{j\ge 1}w_j \tfrac12\big(\cos(\pi j(s-u))-\cos(\pi j(s+u))\big)}\\
&\le \frac1h\Big(w_0+\sum_{j\ge 1}w_j \big(\cos(\pi j(s-u))+\cos(\pi j(s+u))\big)\Big)\\
&=\sum_{j\ge 0}w_j \Psi_j(s)\Psi_j(u).
\end{align*}
By definition and orthonormality of $(\Psi_j)$, we have
\begin{align*}
v_2^2
 &=\sup_{\norm{F}\le 1} \int_0^1\Big\|\Big(\int_u^1\sum_{j\ge 1}w_j \Phi_j(s)\Phi_j(u)\big(\scapro{\Sigma^{1/2}_{>r}(s)}{F(s)}_{HS}I_{d-r}+\Sigma_{>r}^{1/2}(s)F(s)^\top
 \big)\, ds\Big) \\
 &\qquad\qquad\qquad \Sigma_{>r}^{1/2}(u)\Big\|_{HS}^2du\\
&\le \sup_{\norm{F}\le 1} \int_0^h \Big(\int_0^h\Big(\sum_{j\ge 0}w_j \Psi_j(s)\Psi_j(u)\Big) \big(\norm{\Sigma_{>r}^{1/2}(s)}_{HS}+\norm{\Sigma_{>r}^{1/2}(s)}\big)\norm{F(s)}_{HS}\, ds \Big)^2\\
&\qquad \times\norm{\Sigma_{>r}^{1/2}(u)}_{HS}^2\,du\\
&\le  w_0^2 \Big(\norm{\norm{\Sigma_{>r}^{1/2}(\cdot)}_{HS}}_{L^\infty([0,h])}
+\norm{\Sigma_{>r}}_{L^\infty([0,h])}^{1/2}
\Big)^2  \norm{\norm{\Sigma_{>r}^{1/2}(\cdot)}_{HS}^2}_{L^\infty([0,h])}\\
&\le w_0^2(r^{1/2}+1)^2r\norm{\Sigma_{>r}}_{L^\infty([0,h])}^2
\end{align*}
because of $\rank(\Sigma_{>r}(s))\le r$ under ${\mathcal H}_0$.
We thus obtain
\[ \max(v_1, v_2)\le \Big(\norm{(w_j)}_{\ell^2}\vee (r^{1/2}+1)w_0\Big)\sqrt{r+1}\norm{\Sigma_{>r}}_{L^\infty([0,h])}.
\]

For $\sigma^2$ the orthonormality of $(\Phi_j)$ yields with $v\in\R^{d-r}$
\begin{align*}
\sigma^2 &= \sup_{\norm{f}=1} \bnorm{\sum_{j\ge 1} w_j \Big(\int_0^1  \Sigma(t)_{>r}^{1/2}\Phi_j(t)f(t)\,dt\Big)^{\otimes 2}}\\
&= \sup_{\norm{f}=\norm{v}=1} \sum_{j\ge 1} w_j  \Big(\int_0^h\scapro{\Sigma(t)_{>r}^{1/2}f(t)}{v}\Phi_j(t)\,dt\Big)^2\\
&\le \sup_{j\ge 1} w_j\sup_{\norm{f}=\norm{v}=1} \norm{\scapro{\Sigma_{>r}^{1/2}f}{v}}_{L^2([0,h])}^2\\
&\le w_1\norm{\Sigma_{>r}}_{L^\infty([0,h])}.
\end{align*}

Theorem \ref{ThmMatrixDev} with $\delta=1$ and $\lambda_{max}( \Sigma_{>r}^{t,h} )\le \norm{\Sigma_{>r}}_{L^\infty([0,h])}$ thus gives with probability at least $1-\alpha$
\begin{align*}
\lambda_{max}(\tilde\Sigma_{>r}^{t,h})& \le \Big(1+
 \log(\sqrt{7e/2}(d-r))\big(\norm{(w_j)}_{\ell^2}\vee (r^{1/2}+1)w_0\big)\sqrt{r+1} \\
&\qquad  +2w_1\log(\alpha^{-1}) \Big)2\norm{\Sigma_{>r}}_{L^\infty([0,h])}.
\end{align*}
Under ${\mathcal H}_0$ we  now insert $\norm{\Sigma_{>r}}_{L^\infty(I_{t,h})}\le\BIAS_0^\Sigma(h)$ from \eqref{EqBiasRW}
as well as the weight bounds $\norm{(w_j)}_{\ell^2}\le M^{-1/2}$, $w_0\vee w_1\le (4/\pi)M^{-1}$
and we obtain the result after simplification of numerical constants, using $r\le M$.
\end{proof}

\begin{corollary}\label{CorChat}
We have with probability at least $1-\alpha$  under ${\mathcal H}_0$
\begin{align*}
\lambda_{r+1}(\widehat C^{t,h})\le \lambda_{max}(\widehat C_{>r}^{t,h})&
  \le \big(4\BIAS_0^\Sigma(h)+\eps_{n,h}^2\tfrac{M^2}2\big)\times\\
  &\quad \Big(1+
 \log(4(d-r))\Psi(6dM^{-1})
 +21M^{-1}\log(2\alpha^{-1}) \Big).
\end{align*}
\end{corollary}

\begin{proof}
For any $v\in\R^d$ we bound
\begin{align*} \scapro{\widehat C_{>r}^{t,h}v}{v}&=\sum_{j\ge 1} w_j\bscapro{\int_{I_{t,h}}\Phi_j(s)\Sigma_{>r}^{1/2}(s)\,dB(s)+\int_{I_{t,h}}\phi_j(s)dW_{>r}(s)}{v}^2\\
&\le 2\scapro{\tilde\Sigma_{>r}^{t,h}v}{v}+2\scapro{\tilde V_{>r}^{t,h}v}{v},
\end{align*}
using $(A+B)^2\le 2A^2+2B^2$ and inserting the definitions of $\tilde\Sigma_{>r}^{t,h}$, $\tilde V_{>r}^{t,h}$. By \eqref{EqCIL} this shows
\[ \lambda_{r+1}(\widehat C^{t,h})\le \lambda_{max}(\widehat C_{>r}^{t,h})\le 2\lambda_{max}(\tilde\Sigma_{>r}^{t,h})+2\lambda_{max}(\tilde V_{>r}^{t,h}).\]
Hence, applying Propositions \ref{PropVtilde} and \ref{PropSigmatilde} with $\alpha/2$, the result follows by a union bound and a simplification with $(d-r+1)\vee r\le d$.
\end{proof}

\subsection{Lower eigenvalue deviation bounds}\label{secMatrixLB}

The analysis of the test's power under the alternative relies on lower bounds for minimal singular values of Gaussian matrices and relatively abstract entropy and compactness arguments.

\begin{theorem}\label{ThmBernsteinLower}
Let $Y\in\R^{d\times J}$ be a centred Gaussian matrix with $\vek(Y)\sim N(0,A)$  for $A\in\R_{spd}^{dJ\times dJ}$.
Then  we have for $J\ge d$ the lower tail bound
\begin{align*}
\PP\Big(\lambda_{min}\big(YY^\top\big)\le \tau(J-d+2)\lambda_{min}(A)\Big) & \le (2^{d+1}e)^{1/2} \big(2ed \tau\big)^{(J-d+1)/2},\quad \tau\ge 0.
\end{align*}
In particular, for $\delta\in(0,1)$ we have with probability at least $1-\delta$:
\begin{align*}
\lambda_{min}\big(YY^\top\big)\ge e^{-2}2^{-d-2}\delta^{2/(J-d+1)} (J-d+2)\lambda_{min}(A).
\end{align*}
\end{theorem}

\begin{proof}
Introduce the elementary matrices $E_{i,j}:=e_ie_j^\top$ and
$F_{\pm 1}:=(I_{J}-E_{d,d})\pm E_{d,d}\in \R^{J\times J}$ so that $F_1$ is the identity and $F_{-1}$ flips the sign of the $d$th coordinate. Then by the invariance $O\stackrel{d}{=} F_{\pm 1} O$ for $O\sim\HH_J$ with the Haar measure $\HH_J$ on the orthogonal group $O(J)\subset\R^{J\times J}$  we can introduce a random sign flip $F_\eps$ with a Rademacher random variable $\eps$ in the density formula of Lemma S.5 in \cite{RW23}  to obtain
\begin{align*}
&f_A^\Lambda(\lambda_1,\ldots,\lambda_{d})=c\det(A)^{-1/2}\prod_{j=1}^{d}\lambda_j^{(J-d-1)/2}\prod_{i<j}(\lambda_i-\lambda_j)\times\\
&
\quad\iint \E_\eps\Big[\exp\Big(-\frac12\scapro{A^{-1}\vek(O' LF_\eps O)}{\vek(O' LF_\eps O)}_{\R^{d J}} \Big)\Big]\,\HH_J(dO)\HH_d(dO'),
\end{align*}
where $L=(\diag(\lambda_1^{1/2},\ldots,\lambda_d^{1/2}),{\bf 0}_{d\times(J-d)})$.
For $A=A^{(0)}:=\lambda_{min}(A)I_{dJ}$ we use $OO^\top=I_{J}=F_\eps^2$ to obtain
\begin{align*}
&\scapro{(A^{(0)})^{-1}\vek(O' LF_\eps O)}{\vek(O' LF_\eps O)}_{\R^{d J}}\\
&=\lambda_{min}(A)^{-1}\trace\big(O' LF_\eps O O^\top F_\eps L^\top (O')^\top\big)=\lambda_{min}(A)^{-1}
\trace\big( L L^\top \big),
\end{align*}
independently of $\eps,O,O'$. The likelihood ratio $\frac{f_A^\Lambda(\lambda_1,\ldots,\lambda_{d})}{f_{A^{(0)}}^\Lambda(\lambda_1,\ldots,\lambda_{d})}$ between the laws under $A$ and $A^{(0)}$ can then be written as
\[\tilde c \iint\E_\eps\Big[ \exp\Big(\tfrac12\scapro{((A^{(0)})^{-1}-A^{-1})\vek(O' LF_\eps O)}{\vek(O' LF_\eps O)}_{\R^{dJ}}\Big)\Big]\,\HH_J(dO)\HH_d(dO')
\]
with some constant $\tilde c>0$. Noting
$LF_\eps=(\sum_{i\le d-1}\lambda_i^{1/2}E_{i,i})+\eps \lambda_{d}^{1/2}E_{d,d}$,
the exponent is a quadratic form $Q(\eps\lambda_{d}^{1/2})$ in $\eps \lambda_{d}^{1/2}$. Since $(A^{(0)})^{-1}-A^{-1}$ is positive semi-definite by definition, we can write $Q(x)=a+bx+cx^2$ with some $a,c\ge 0$ and $b\in\R$. This shows that
\[ \E_\eps[\exp(Q(\eps \lambda_{d}^{1/2}))]= \exp(a+c\lambda_{d})\cosh(b\lambda_{d}^{1/2})\]
is always increasing in $\lambda_{d}\ge 0$. Hence,  the likelihood ratio $f_A^\Lambda/f_{A^{(0)}}^\Lambda$ is increasing in $\lambda_{d}$. Integrating $\lambda_i$ for $i\le d-1$ out, we thus get the stochastic order
\[ \PP\big(\lambda_{min}(YY^\top)\le \tau\big)\le \PP^{(0)}\big(\lambda_{min}(YY^\top)\le \tau\big), \quad \tau\ge 0,
\]
where $\PP^{(0)}$ denotes the law under $A^{(0)}$. Writing $Y=\lambda_{min}(A)^{1/2}(\zeta_1,\ldots,\zeta_J)\in\R^{d\times J}$ with $\zeta_j\sim N(0,I_d)$ i.i.d. under $\PP^{(0)}$, we conclude
\begin{align*}
\PP\Big(\lambda_{min}(YY^\top)\le \lambda_{min}(A)\tau\Big)
&\le\PP\Big(\lambda_{min}\Big(\sum_{j=1}^{J}\zeta_j^{\otimes 2}\Big)\le \tau\Big), \quad \tau\ge 0.
\end{align*}

Using the bound for the density of $\lambda_{min}$ in the proof of Prop. 5.1 in \citet{Edelman1988}, we obtain for $\tau\ge 0$
\begin{align*}
&\PP\Big(\lambda_{min}\Big(\sum_{j=1}^{J}\zeta_j^{\otimes 2}\Big)\le \tau\Big)\\
&\le \frac{ 2^{-(J-d+1)/2} \Gamma(1/2)\Gamma((J+1)/2)} {\Gamma(d/2)\Gamma((J-d+1)/2)\Gamma((J-d+2)/2)} \int_0^\tau \lambda^{(J-d-1)/2}e^{-\lambda/2}d\lambda\\
&\le \frac{2^{-(J-d+1)/2} \Gamma(1/2)\Gamma((J+1)/2)} {\Gamma(d/2)\Gamma((J-d+3)/2)\Gamma((J-d+2)/2)} \tau^{(J-d+1)/2}.
\end{align*}
By Euler's product formula $\Gamma(z)=z^{-1}\prod_{n\ge 1}(1+1/n)^z(1+z/n)^{-1}$ we can derive for all $0<y<x$ the bound
\[ \frac{\Gamma(x)}{\Gamma(y)\Gamma(x-y)}\le \frac{\Gamma(x)}{\Gamma(x/2)^2}=\frac x4 \prod_{n\ge 1}\Big(1+\Big(4(n/x)^2+4n/x\Big)^{-1}\Big)\le \frac x4 2^x,\]
where the last inequality follows by taking logarithms and by comparison with $\int_0^\infty\log(1+(4x^2+4x)^{-1})\,dx=\log 2$. With $x=(J+1)/2$ and $y=d/2$ it follows that
\[ \PP\Big(\lambda_{min}\Big(\sum_{j=1}^{J}\zeta_j^{\otimes 2}\Big)\le \tau\Big) \le \frac{2^{d/2} \Gamma(1/2)} {((J-d+1)/2)\Gamma((J-d+2)/2)} \tau^{(J-d+1)/2}.
\]
Now using the Stirling formula bound $\Gamma(z)\ge \sqrt{2\pi} z^{z-1/2}e^{-z}$ for $z>0$ \citep{artin2015} we arrive at
\[ \PP\Big(\lambda_{min}\Big(\sum_{j=1}^{J}\zeta_j^{\otimes 2}\Big)\le \tau\Big) \le \frac{ (2^{d+1}e)^{1/2}} {J-d+1} \Big(\frac{2e\tau}{J-d+2}\Big)^{(J-d+1)/2},\quad \tau\ge 0.
\]
Replacing $\tau$ by $(J-d+2)\tau$ and some direct simplifications yield the assertion.
\end{proof}

\begin{proposition}\label{PropEigenspace}
For $J\in\N$ there is a universal constant $\underline{C}_J>0$ (depending otherwise only on $d,r$) such that for any interval $I_{t,h}$ and $\Sigma$ with $\lambda_{r+1}(\Sigma(s))\ge\underline\lambda_{r+1}$, $s\in I_{t,h}$, there is an $(r+1)$-dimensional subspace $V_{r+1}\subset\R^d$ such that for all $v_j\in\R^{r+1}$, $j=1,\ldots,J$:
\begin{align*}
&\int_{I_{t,h}} \lambda_{min}\big(P_{V_{r+1}}\Sigma(s)|_{V_{r+1}}\big)  \Big\|\sum_{j=1}^J w_{j}^{1/2}\Phi_{j}(s)v_{j}\Big\|^2ds
 \ge \underline{C}_J\tfrac14\underline\lambda_{r+1}\Big(\min_{j=1,\ldots,J} w_j\Big)\sum_{j=1}^J\norm{v_j}^2.
\end{align*}

\end{proposition}

\begin{proof}
Without loss of generality consider $I_{t,h}=[0,h]$ with $\Phi_j(s)=\sqrt{2/h}\sin(j\pi s/h)$.
From $\lambda_{r+1}(\Sigma(s))\ge\underline\lambda_{r+1}$ we infer by the variational characterisation of eigenvalues that there is an $(r+1)$-dimensional subspace $V_{r+1}(s)$ such that $\lambda_{min}(P_{V_{r+1}(s)}\Sigma(s)|_{V_{r+1}(s)})\ge\underline\lambda_{r+1}$.
By \citet{Pajor1998} there is a finite set $\cal V$ of $(r+1)$-dimensional subspaces of $\R^d$ with $\min_{V\in{\cal V}}\norm{P_{V_{r+1}}-P_V}\le 1/2$ for all $(r+1)$-dimensional subspaces of $\R^d$. By a simple counting argument, compare the proof of Cor. S.9 in \cite{RW23}, there is $V_{r+1}\in{\cal V}$ and a Borel set $B_{r+1}\subset [0,h]$ with $\lambda(B_{r+1})\ge h/\abs{\cal V}=:ch$  ($c\in(0,1)$ depends only on $r$ and $d$) and
\[\forall s\in B_{r+1}:\,\lambda_{min}(P_{V_{r+1}(s)}\Sigma(s)|_{V_{r+1}})\ge\tfrac14\underline\lambda_{r+1}.\]

Next, note that
\begin{align*}
T_{v,J}(s)&:=\Big\|\sum_{j=1}^J w_{j}^{1/2}\Phi_{j}(s)v_{j}\Big\|^2=\frac2h \sum_{j,j'=1}^J w_{j}^{1/2}w_{j'}^{1/2}\scapro{v_j}{v_{j'}}\sin(j\pi s/h)\sin(j'\pi s/h)
\end{align*}
is a nonnegative trigonometric polynomial of degree $2J$.

For the space $\T_{2J}$ of trigonometric polynomials $\sum_{j=-2J}^{2J}\alpha_j e^{ijs}$ of degree at most $2J$ consider the linear functional
\[ \Xi:\T_{2J}\times L^1([0,1])\to\R;\quad \Xi(T,f)=\int_0^1T(s)f(s)\,ds.
\]
Introduce (with $c\in(0,1)$ from above) the  subsets
\begin{align*}
  A &=\Big\{T\in\T_{2J}\,\Big|\,T\ge 0,\,\int_0^1T=1\Big\},\,
  F=\Big\{f\in L^1([0,1])\,\Big|\, 0\le f\le 1,\,\int_0^1f\ge c\Big\}.
\end{align*}
It is well known that any non-zero trigonometric polynomial $t\in\T_N$ of degree $N$ has at most $2N$ zeroes. The nice argument is that $t(s)=0$ implies $p(e^{is})=0$ for the corresponding complex polynomial $p(z)=z^N\sum_{j=-N}^N\alpha_j z^j$, which by the fundamental theorem of algebra has at most $2N$ zeroes in $\C$. Hence, $T>0$ holds Lebesgue-almost everywhere for $T\in\T_{2J}$, implying $\Xi(T,f)>0$ for any $T\in\T_{2J}$, $f\in F$. We shall employ a compactness argument to show
\begin{equation}\label{EqCJ}
\underline C_J:=\inf_{f\in F,T\in A}\int_0^1 T(s)f(s)\,ds>0.
\end{equation}
We equip the finite-dimensional space $\T_{2J}$ with the $L^\infty$-norm. Then $\Xi$ is obviously norm continuous and thus also continuous with respect to the weak topology on $\T_{2J}\times L^1([0,1])$ by the definition of weak convergence.
The subset $A\subset \T_{2J}$ is obviously closed in $L^\infty([0,1])$-norm and bounded in $L^1$-norm. By the equivalence of norms on finite-dimensional spaces, $A$ is  compact.  The family $F$ of functions in $L^1([0,1])$ is uniformly integrable (even uniformly bounded), hence weakly relatively compact by the Dunford-Pettis criterion \cite[Thm. VIII.6.9]{werner2007}. Moreover, $F$ is weakly closed as intersection of weakly closed half spaces:
\[ F=\bigcap_{\stackrel{g\in L^\infty([0,1])}{0\le g\le 1}}\Big\{f\in L^1([0,1])\,\Big|\, \int_0^1 fg\in [0,1]\Big\}\bigcap \Big\{f\in L^1([0,1])\,\Big|\, \int_0^1 f\ge c\Big\}.
\]
Hence, $A\times F\subset \T_{2J}\times L^1([0,1])$ is compact as a product of compact spaces and the infimum in \eqref{EqCJ} is attained, whence positive.

Applying \eqref{EqCJ} with indicators $f={\bf 1}_{B'}$, we conclude
\begin{align*}
& \inf_{B\subset [0,h],\lambda(B)\ge ch}\int_B  \Big\|\sum_{j=1}^J w_{j}^{1/2}\Phi_{j}(s)v_{j}\Big\|^2\,ds
=\inf_{B'\subset [0,1],\lambda(B')\ge c}\int_{B'}   T_{v,J}(hu)h\,du\\
&\ge \underline C_J\int_0^1 T_{v,J}(hu)h\,du=\underline C_J \sum_{j=1}^J w_j\norm{v_j}^2.
\end{align*}
Using $\lambda_{r+1}(\Sigma(s))\ge\underline\lambda_{r+1}$, $s\in [0,h]$, and  $V_{r+1}$, $B_{r+1}$ from above,  we  bound
\begin{align*}
&\int_0^h \lambda_{min}(P_{V_{r+1}}\Sigma(s)|_{V_{r+1}})  \Big\|\sum_{j=1}^J w_{j}^{1/2}\Phi_{j}(s)v_{j}\Big\|^2\,ds\\
&\ge \int_{B_{r+1}} \frac{\underline\lambda_{r+1}}{4}  \Big\|\sum_{j=1}^J w_{j}^{1/2}\Phi_{j}(s)v_{j}\Big\|^2\,ds\ge \tfrac14 \underline{C}_J\underline\lambda_{r+1}\Big(\min_{j=1,\ldots,J} w_j\Big)\sum_{j=1}^J\norm{v_j}^2,
\end{align*}
 and the result follows.
\end{proof}

\begin{corollary}\label{CorLambdaLower}
Suppose $\lambda_{r+1}(\Sigma(s))\ge \underline\lambda_{r+1}>0$ for $s\in I_{t,h}$. Then
\begin{equation}\label{EqLambdaLower} \PP\Big(\lambda_{r+1}\big(\widehat C^{t,h}\big)\le \tau M^{-1}(\underline\lambda_{r+1}+\eps_{n,h}^2) \Big) \le C_J  \tau^{(J-r)/2},\quad \tau\ge 0,
\end{equation}
holds for $r+1\le J\le M$ with $C_J>0$ only depending on $J$, $d$ and $r$.

In particular, we have for $R\to\infty$
\[\lim_{R\to\infty}\inf_{\Sigma\in {\mathcal H}_1(I_{t,h},RM\kappa_{\alpha,i})} \PP_{\Sigma}(\phi_{\alpha,i}=1)=1.\]
\end{corollary}

\begin{remark}\label{RemH1beta}
Ideally, we would like to have \eqref{EqLambdaLower} without the factor $M^{-1}$ and for $J=M$ with $C_M$ satisfying $\inf_MC_M>0$. The inherent problem is indeed to find a common $(r+1)$-dimensional subspace such that $\Sigma(s)$ for all $s\in I_{t,h}$ has minimal eigenvalue of order larger than $\underline\lambda_{r+1}$ on this subspace, for which we use Proposition \ref{PropEigenspace} in the proof. For a rank $r=d-1$ this is trivially the full space $\R^{d}$ and we get rid of the factor $M^{-1}$. Also, under the additional regularity assumption $\Sigma\in C^\beta(L)$ this is feasible for $\underline\lambda_{r+1}\ge 2Lh^\beta$, using $\norm{\Sigma(s)-\Sigma(u)}\le Lh^\beta$ for $s,u\in I_{t,h}$. This  covers the critical values without spectral gap, but not those with spectral gap where typically $\underline\lambda_{r+1}\thicksim h^{2\beta}$. In the latter case we cannot expect to find a suitable common subspace given the possibility of several indices $k\ge r+1$ with eigenvalues $\lambda_k(\Sigma(s))\approx \lambda_{r+1}(\Sigma(s))$, compare the discussion of the spiked covariance model in \citet{reiss2020}.
\end{remark}

\begin{proof}
We consider $Y=(w_1^{1/2}P_{V_{r+1}}S_{1},\ldots,w_J^{1/2}P_{V_{r+1}}S_{J})\in \R^{(r+1)\times J}$ in Theorem \ref{ThmBernsteinLower} with $S_j$ from \eqref{EqSj} and the projection $P_{V_{r+1}}$ onto the common $(r+1)$-dimensional subspace $V_{r+1}$, provided by Proposition \ref{PropEigenspace}. Then taking duality with respect to the Hilbert-Schmidt scalar product, the covariance matrix $A\in\R^{(r+1)J\times (r+1)J}$ of $\vek(Y)$ satisfies for $H=(H_1,\ldots,H_J)\in \R^{(r+1)\times J}$
\[ \lambda_{min}(A)=\inf_{\norm{H}_{HS}=1} \E[\scapro{Y}{H}_{HS}^2]= \inf_{\sum_{j=1}^J\norm{H_j}^2=1} \E\Big[\Big(\sum_{j=1}^Jw_j^{1/2}\scapro{P_{V_{r+1}}S_{j}}{H_j}\Big)^2\Big].
\]
Writing $\Sigma(s)_{V_{r+1}}=P_{V_{r+1}}\Sigma(s)|_{V_{r+1}}$, we can lower bound
\begin{align*}
 & \E\Big[\Big(\sum_{j=1}^Jw_j^{1/2}\scapro{P_{V_{r+1}}S_{j}}{H_j}\Big)^2\Big]\\
&= \sum_{j,j'=1}^Jw_j^{1/2}w_{j'}^{1/2} H_j^\top P_{V_{r+1}}\E\big[S_{j}S_{j'}^\top \big] P_{V_{r+1}}^\top H_{j'}^\top\\
&= \sum_{j,j'=1}^J  \int_{I_{t,h}} \scapro{w_j^{1/2}\Phi_j(s)\Sigma(s)_{V_{r+1}}H_j}{w_{j'}^{1/2}\Phi_{j'}(s)H_{j'}}\,ds
 +\sum_{j=1}^J w_j j^2\eps_{n,h}^2\norm{H_j}^2\\
&\ge \int_{I_{t,h}} \lambda_{min}(\Sigma(s)_{V_{r+1}})\Big\|\sum_{j=1}^J w_{j}^{1/2}\Phi_{j}(s)H_{j}\Big\|^2ds
 +\eps_{n,h}^2\sum_{j=1}^J w_j j^2\norm{H_j}^2.
\end{align*}
Consequently, from Proposition \ref{PropEigenspace} we deduce
\[  \lambda_{min}(A)\ge \frac{\underline{C}_J}{4}\Big(\min_{j=1,\ldots,J}w_j\Big)\underline\lambda_{r+1}+\eps_{n,h}^2\min_{j=1,\ldots,J} (w_j j^2).
\]
For $J\le M$ we have $w_j\thicksim M^{-1}$ and thus
\[  \lambda_{min}(A)=\inf_{\norm{H}_{HS}=1} \E[\scapro{Y}{H}_{HS}^2]\gtrsim \underline{C}_J M^{-1}(\underline\lambda_{r+1}+\eps_{n,h}^2).\]
 Hence, Theorem \ref{ThmBernsteinLower} yields for $M\ge J\ge r+1$ with a  constant $C_J>0$ (only depending on $J$, $d$ and $r$)
 \begin{align*}
\PP\Big(\lambda_{min}\big(YY^\top\big)\le  \tau M^{-1}(\underline\lambda_{r+1}+\eps_{n,h}^2)\Big) & \le C_J  \tau^{(J-r)/2}.
\end{align*}
By the variational characterisation of eigenspaces  this yields
  \begin{align*}
&\PP\Big(\lambda_{r+1}\big(\widehat C^{t,h}\big)\le \tau M^{-1}(\underline\lambda_{r+1}+\eps_{n,h}^2) \Big) \\
&\le \PP\Big(\lambda_{r+1}\Big(\sum_{j=1}^J w_jS_{j}^{\otimes 2}\Big)\le  \tau M^{-1}(\underline\lambda_{r+1}+\eps_{n,h}^2)\Big)
  \le C_J  \tau^{(J-r)/2}
\end{align*}
for any $r+1\le J\le M$, which proves the first statement.

Now with $\underline\lambda_{r+1}=RM\kappa_{\alpha,i}$, we deduce in particular for $\tau=R^{-1}$ and $J=r+1$
\[ \PP\big(\lambda_{r+1}\big(\widehat C^{t,h}\big)\le \kappa_{\alpha,i} \big) \le C_{r+1}  R^{-1/2}\to 0\]
 as $R\to\infty$, uniformly over $\Sigma\in{\mathcal H}_1(I_{t,h},RM\kappa_{\alpha,i})$.
\end{proof}

\subsection{Proof of the lower bound}\label{SecProofLB}

\begin{proof}[Proof of Theorem \ref{ThmLowerBoundBlock}]
We  treat without loss of generality the case $\eta=1$, the general case follows via replacing $n$ by $\eta^{-2}n$.
Moreover, we shall only consider $\underline\lambda_r\ge (L^2 n^{-\beta})^{1/(\beta+2)}$ because for $\underline\lambda_r\le (L^2 n^{-\beta})^{1/(\beta+2)}$ the detection rate is always $v_n= L^{2/(\beta+2)}n^{-\beta/(\beta+2)}$ and the null hypotheses ${\mathcal H}_0$ are nested with respect to the spectral gap condition.

First, let $d=2$, $r=1$ and omit the index $n$ in all quantities. For $h\in(0,(2t)\wedge \frac27(1-t)\wedge (\underline\lambda_r/L)^{1/\beta}]$, to be fixed below,   we consider the two alternatives (inspired by Example 2.2 in \citet{RW23})
\begin{align*}
\Sigma_0(s) &=\begin{pmatrix} \underline\lambda_r & L h^\beta f((s-t)/4h)\\ L h^\beta f((s-t)/4h) & \underline\lambda_r^{-1}L^2h^{2\beta}f^2((s-t)/4h)\end{pmatrix},\\
 \Sigma_1(s)&=\begin{pmatrix} \underline\lambda_r & 0\\  0& \underline\lambda_r^{-1}L^2h^{2\beta}f^2((s-t)/4h)\end{pmatrix},
\end{align*}
satisfying $\Sigma_0=\Sigma_1^{1/2}\begin{pmatrix} 1 & 1\\ 1&1\end{pmatrix}\Sigma_1^{1/2}$, with
\[ f(u)=(4\pi)^{-1}\cos\big(2\pi (u-1/8)\big){\bf 1}(u\in[-\tfrac18,\tfrac78]).\]
Then $f$ and $f^2$ are $\beta$-H\"older continuous with constant $1/2$ and by $Lh^{\beta}\le\underline\lambda_r$ this yields $\Sigma_0,\Sigma_1\in C^\beta([0,1],L)$. Moreover, the ordered eigenvalues are
\begin{align*}
\lambda_1(\Sigma_0(s))&=\underline\lambda_r+\underline\lambda_r^{-1}L^2h^{2\beta}f^2((s-t)/4h)
\ge\underline\lambda_r,\qquad \lambda_2(\Sigma_0(s))=0,\\
\lambda_1(\Sigma_1(s))&=\underline\lambda_r,\qquad \lambda_2(\Sigma_1(s))=\underline\lambda_r^{-1}L^2h^{2\beta}f^2((s-t)/4h).
\end{align*}
 We conclude
\begin{equation}\label{EqLBclasses}
\Sigma_0\in {\mathcal H}_0([0,1],\beta,L,\underline\lambda_r),\quad \Sigma_1\in {\mathcal H}_1(I_{t,h},(32\pi^2)^{-1}\underline\lambda_r^{-1}L^2h^{2\beta})\cap C^\beta([0,1],L)
\end{equation}
due to $f^2(u)\ge\frac1{32\pi^2}$ for $u\in[0,1/4]$.
Below it will be important that $\int_0^{s}\Sigma_0(u)du=\int_0^{s}\Sigma_1(u)du$ holds for $s\le \underline t_h:=t-h/2$ and $s\ge \bar t_h:=t+7h/2$ due to $\int f(u)du=0$.

We observe \eqref{EqY} for $\Sigma\in\{\Sigma_0,\Sigma_1\}$.
Since $\Sigma_0(u)=\Sigma_1(u)$ holds for $u$ outside $[\underline t_h,\bar t_h]$  and $X(t)$ is Markovian, the observation laws for $\Sigma\in\{\Sigma_0,\Sigma_1\}$ coincide on $[0,\underline t_h]\cup[\bar t_h,1]$ conditional on the boundary values $X(\underline t_h)$, $X(\bar t_h)$. By adding the observation of the values $X(\underline t_h)$, $X(\bar t_h)$ the model becomes more informative and the lower bound can only become smaller. Note that $(X(\underline t_h),X(\bar t_h))$ is a centred Gaussian random vector with  covariance independent of $i$ in $\Sigma=\Sigma_i$.
 Hence, it suffices to prove the lower bound for the models of observing for $s\in [\underline t_h,\bar t_h]$
\begin{align*}
  d\tilde Y(s)=&dY(s)-\big(\tfrac{\bar t_h-s}{\bar t_h-\underline t_h} X(\underline t_h)+\tfrac{s-\underline t_h}{\bar t_h-\underline t_h} X(\bar t_h)\big)ds=\tilde X(s)ds+n^{-1/2}dW(s) \\
&\text{with }\tilde X(s)=X(s)-\tfrac{\bar t_h-s}{\bar t_h-\underline t_h} X(\underline t_h)-\tfrac{s-\underline t_h}{\bar t_h-\underline t_h}X(\bar t_h).
\end{align*}
The bridge process $\tilde X$ is centred Gaussian with covariance function $c_i$ for $u,v\in[\underline t_h, \bar t_h]$ under $\Sigma=\Sigma_i$, $i=0,1$:
\begin{align*}
&c_i(u,v) =A_i(u\wedge v)-\tfrac{u-\underline t_h}{\bar t_h-\underline t_h} A_i(v)-\tfrac{v-\underline t_h}{\bar t_h-\underline t_h} A_i(u)+\tfrac{u-\underline t_h}{\bar t_h-\underline t_h}\tfrac{v-\underline t_h}{\bar t_h-\underline t_h}A_i(\bar t_h) \\
&\text{where }A_i(s) =\E_{\Sigma_i}[(X(s)-X(\underline t_h))^{\otimes 2} ] =\int_{\underline t_h}^{s}\Sigma_i(u)\,du.
\end{align*}
Let $Q_i: L^2([\underline t_h,\bar t_h];\R^2)\to L^2([\underline t_h,\bar t_h];\R^2)$ denote the covariance operator of the cylindrical Gaussian measure on $L^2([\underline t_h,\bar t_h];\R^2)$ generated by observing $d\tilde Y$ with $\Sigma=\Sigma_i$.
Then by independence of $B$ and $W$ we have $Q_i=C_i+\frac1n \Id$ with the identity operator $\Id$ and $C_i$ the covariance operator of $\tilde X$ described by the covariance function $c_i$.

We write any linear operator $T$ on $L^2([\underline t_h,\bar t_h];\R^2)$ in block-matrix notation
\[T=\begin{pmatrix} T_{11}& T_{12}\\ T_{21} & T_{22}\end{pmatrix} \text{ with } T\begin{pmatrix} f_1\\ f_2\end{pmatrix}=\begin{pmatrix} T_{11}f_1+T_{12}f_2\\ T_{21}f_1+T_{22}f_2\end{pmatrix}
\]
for $T_{ij}:L^2([\underline t_h,\bar t_h])\to L^2([\underline t_h,\bar t_h])$.
By partial integration we obtain the well known representation
\[ C_i=\begin{pmatrix} J^\ast & 0\\ 0 & J^\ast \end{pmatrix}M_{\Sigma_i}\begin{pmatrix} J & 0\\ 0 & J \end{pmatrix}\text{ with } Jf(s):=\int_{\underline t_h}^{\bar t_h} \big({\bf 1}(u\ge s)-\tfrac{u-\underline t_h}{\bar t_h-\underline t_h}\big)f(u)\,du
\]
and multiplication operators $M_{\Sigma_i} f(s):=\Sigma_i(s)f(s)$, $s\in [\underline t_h,\bar t_h]$, satisfying $M_{\Sigma_0}=M_{\Sigma_1^{1/2}}\begin{pmatrix} \Id & \Id\\\Id &\Id\end{pmatrix}M_{\Sigma_1^{1/2}}$.
Consequently, we find
\[C_1=D^\ast D, \, C_0=D^\ast\begin{pmatrix} \Id &\Id\\ \Id &\Id\end{pmatrix}D\text{ with }
D=\begin{pmatrix} M_{\underline\lambda_r} J&0\\ 0& M_{\underline\lambda_r^{-1}L^2h^{2\beta}f^2((\cdot-t)/8h)}J\end{pmatrix}.
\]

As for the scalar case \citep{reiss2011} the squared Hellinger distance can be bounded by a Hilbert-Schmidt norm and we obtain
\begin{align}
 H^2(N(0,Q_0),N(0,Q_1)) &\le \norm{Q_1^{-1/2}(Q_1-Q_0)Q_1^{-1/2}}_{HS}^2\nonumber\\
 &=\bnorm{(C_1+\tfrac1n \Id)^{-1/2}D^\ast\begin{pmatrix} 0 & \Id\\ \Id & 0\end{pmatrix}D(C_1+\tfrac1n \Id)^{-1/2}}_{HS}^2\nonumber\\
 &=2\norm{(C_{1,11}+\tfrac1n \Id)^{-1/2}D_{11}^\ast D_{22}(C_{1,22}+\tfrac1n \Id)^{-1/2}}_{HS}^2\nonumber\\
 &\le 2\norm{D_{22}n^{1/2}}_{HS}^2=2n\trace(C_{1,22}),\label{EqHellinger}
 \end{align}
where in the last inequality $\norm{D_{11}(C_{1,11}+\tfrac1n\Id)^{-1} D_{11}^\ast}\le 1$ and $\norm{(C_{1,22}+\tfrac1n \Id)^{-1}}\le n$ were used.

For positive trace-class operators on $L^2$ the trace can be calculated explicitly via the kernel (see e.g. Example VI.5(d) in \citet{werner2007}):
\begin{align}
\trace(C_{1,22}) &= \int_{\underline t_h}^{\bar t_h}c_{1,22}(s,s) ds\nonumber\\
 &=\int_{\underline t_h}^{\bar t_h}\Big(\frac{(\bar t_h-s)^2}{(\bar t_h-\underline t_h)^2}\int_{\underline t_h}^s \Sigma_{1,22}(u)du+\frac{(s-\underline t_h)^2}{(\bar t_h-\underline t_h)^2}\int_s^{\bar t_h}\Sigma_{1,22}(u)du\Big)ds\nonumber\\
 &\le (\bar t_h-\underline t_h)^2\norm{\Sigma_{1,22}}_{\infty}=\tfrac{4}{\pi^2}\underline\lambda_r^{-1}L^2h^{2\beta+2}.\label{EqTrace}
\end{align}
We conclude
\[  H^2(N(0,Q_0),N(0,Q_1))\le \underline\lambda_r^{-1}L^2h^{2\beta+2}n.\]
For any test $\phi_n$, based on observing $d\tilde Y$, the sum of error probabilities satisfies with total variation norm $TV$ and Hellinger distance $H$ \citep{tsybakov2008}
\begin{align*}
 \E_{\Sigma_0}[\phi_n]+\E_{\Sigma_1}[1-\phi_n] &=1-TV(N(0,Q_0),N(0,Q_1))\\
 &\ge 1-H(N(0,Q_0),N(0,Q_1))\\
 &\ge 1- \underline\lambda_r^{-1/2}Lh^{\beta+1}n^{1/2}.
\end{align*}
This bound equals $\alpha'$ when
\begin{equation}\label{Eqh}
h=\big((1-\alpha')^2\underline\lambda_{r}L^{-2}n^{-1}\big)^{1/(2\beta+2)}\thicksim \big(\underline\lambda_r^{-1} L^2 n\big)^{-1/(2\beta+2)}.
\end{equation}
Since we are in the case $\underline\lambda_{r}\ge (L^2n^{-\beta})^{1/(\beta+2)}$,  the condition $h\le(\underline\lambda_{r}/L)^{1/\beta}$ is always satisfied. We obtain for the $(r+1)$st (i.e., second) eigenvalue of $\Sigma(s)$, $s\in I_{t,h}$, under the alternative
\[ \lambda_{r+1}(\Sigma_1(s))\ge(32\pi^2)^{-1}\underline\lambda_{r}^{-1}L^2h^{2\beta}= \tfrac{(1-\alpha')^{2\beta/(\beta+1)}}{32\pi^2}(L^2/\underline\lambda_{r})^{1/(\beta+1)}n^{-\beta/(\beta+1)}.
\]
For $c_{\alpha'}=(32\pi^2)^{-1} (1-\alpha')^{2\beta/(\beta+1)}$ we conclude
\begin{equation}\label{EqLB}
\liminf_{n\to\infty}\inf_{\phi_n}\Big(\sup_{\Sigma\in {\mathcal H}_0([0,1],\beta,L,\underline\lambda_{r})}\E_{\Sigma}[\phi_n]+\sup_{\Sigma\in {\mathcal H}_1(I_{t,h},c_{\alpha'} v_n)}\E_{\Sigma}[1-\phi_n]\Big)\ge \alpha',
\end{equation}
 where the infimum is taken over all tests $\phi_n$ based on observing $d\tilde Y$ and a fortiori also for all tests $\phi_n$ based on observing $dY$ in \eqref{EqY}. For $h'\in(0,h)$ we have the inclusion
 ${\mathcal H}_1(I_{t,h},c_{\alpha'} v_n)\subset {\mathcal H}_1(I_{t,h'},c_{\alpha'} v_n)$ so that the left-hand side in \eqref{EqLB} can only be larger for $h'$ instead of $h$, defined in \eqref{Eqh}. Consequently, \eqref{EqLB} holds for any $h\in(0,c_{\alpha'}^{1/2\beta}(\underline\lambda_r^{-1} L^2 n)^{-1/(2\beta+2)}]$, as asserted.

Finally, this lower bound  generalises to $d>2$ and rank $1\le r\le d-1$ by considering covariances $\Sigma_i$ with the above prescription in coordinates $1$ and $r+1$, putting $\underline\lambda_r$ on the diagonal of $\Sigma_i$ for coordinates $2,\ldots,r$ and setting all other entries to zero. This follows immediately from the fact that the coordinates $1$ and $r+1$ of $dY$ form a sufficient statistics for the binary statistical model in $\{\Sigma_0,\Sigma_1\}$.
\end{proof}

\end{appendix}

\bibliographystyle{apalike2}
\bibliography{HF-RankTest}

\end{document}